\documentclass[preprint,12pt]{elsarticle}

\usepackage{graphicx}
\usepackage{amssymb}
\usepackage{amsthm}

\usepackage{amsmath}
\usepackage{hyphenat}

\usepackage{examplep}

\usepackage[normalem]{ulem}
\usepackage{xcolor}
\usepackage{cancel}
\usepackage{tabularx}
\usepackage{float}
\usepackage{xspace}
\usepackage{color}
\usepackage[hidelinks]{hyperref} %hidelinks to get rid of the square around hyperlink
\hypersetup{
	colorlinks}%,
\usepackage[font=scriptsize]{caption} 
\newcolumntype{C}{>{\Centering\arraybackslash}X}
\newcommand{\ii}{\operatorname{i}}
\newcommand{\Comma}{\: ,}
\newcommand{\period}{\: .}
\newcommand{\semicolon}{\: ;}
\newcommand{\Ai}{\operatorname{Ai}}
\newcommand{\Bi}{\operatorname{Bi}}
\newcommand{\MATLAB}{\textsc{Matlab}\xspace}
\newcommand\redsout{\bgroup\markoverwith{\textcolor{red}{\rule[0.5ex]{2pt}{0.4pt}}}\ULon}

\let\oldfrac\frac% Store \frac
\renewcommand{\frac}[2]{%
	\mathchoice
	{\oldfrac{#1}{#2}}% display style
	{\oldfrac{#1}{#2}} % text style
	{\oldfrac{#1}{#2}} % script style
	{#1/#2}% script-script style
}

\newtheorem{Remark}{Remark}[section]

\numberwithin{equation}{section}

\newcounter{bla}

\journal{Journal of Computational and Applied Mathematics}

\begin{document}

\begin{frontmatter}



\title{WKB-based scheme with adaptive step size control for the Schrödinger equation in the highly oscillatory regime}


\author{Jannis Körner}
\ead{jannis.koerner@tuwien.ac.at}
\author{Anton Arnold\corref{author}}
\ead{anton.arnold@tuwien.ac.at}
\author{Kirian Döpfner}
\ead{kirian.doepfner@gmail.com}

\cortext[author] {Corresponding author.}
\address{Institute of Analysis and Scientific Computing, Technische Universität Wien, Wiedner Hauptstr. 8-10, 1040 Wien, Austria}

\begin{abstract}
This paper is concerned with an efficient numerical method for solving the 1D stationary Schrödinger equation in the highly oscillatory regime. Being a hybrid, analytical-numerical approach it does not have to resolve each oscillation, in contrast to standard schemes for ODEs. We build upon the WKB-based (named after the physicists Wentzel, Kramers, Brillouin) marching method from \cite{ABN11} and extend it in two ways: By comparing the $\mathcal{O}(h)$ and $\mathcal{O}(h^{2})$ methods from \cite{ABN11} we design an adaptive step size controller for the WKB method. While this WKB method is very efficient in the highly oscillatory regime, it cannot be used close to turning points. Hence, we introduce for such regions an automated methods switching, choosing between the WKB method for the oscillatory region and a standard Runge-Kutta-Fehlberg 4(5) method in smooth regions.

A similar approach was proposed recently in \cite{HLH16,AHLH20}, however, only for an $\mathcal{O}(h)$-method. Hence, we compare our new strategy to their method on two examples (Airy function on the spatial interval $[0,\,10^{8}]$ with one turning point at $x=0$ and on a parabolic cylinder function having two turning points), and illustrate the advantages of the new approach w.r.t.\ accuracy and efficiency.
\end{abstract}


\end{frontmatter}

%\begin{small}
%\begin{thebibliography}{0}
%\bibitem{1}Reference 1         % This list should only contain those items referenced in the                 
%\bibitem{2}Reference 2         % Program Summary section.   
%\bibitem{3}Reference 3         % Type references in text as [1], [2], etc.
%                               % This list is different from the bibliography at the end of 
%                               % the Long Write-Up.
%\end{thebibliography}
%\end{small}

\section{Introduction}\label{chap:Intro}
This paper deals with the numerical solution of the highly oscillatory 1D Schrödinger equation
\begin{align} \label{schroedinger}
	\varepsilon^{2}\varphi^{\prime\prime}(x) + a(x) \varphi(x) = 0 \Comma\quad x\in\mathbb{R} \period
\end{align}
Here, $0<\varepsilon\ll 1$ is the rescaled Planck constant ($\varepsilon:=\frac{\hbar}{\sqrt{2m}}$) and $\varphi$ the (possibly complex valued) Schrödinger wave function. The real valued coefficient function $a(x)$ is related to the potential. We shall assume here that it is bounded away from zero, i.e.\ $a(x)\geq \tau$ for some $\tau > 0$. The (local) wave length of a solution $\varphi$ to (\ref{schroedinger}) is given by $\lambda(x)=\frac{2\pi\varepsilon}{\sqrt{a(x)}}$. Hence, for small values of $\varepsilon$ the solution is highly oscillatory, especially in the semi-classical limit $\varepsilon \to 0$.

Oscillatory problems like (\ref{schroedinger}) appear in a wide range of applications, e.g., quantum mechanics, electron transport in semiconductor devices, and acoustic scattering.
% Concrete application
For instance, the state of an electron that is injected with the prescribed energy $E$ from the right boundary into an electronic device (e.g., diode), modeled on the interval $[0,1]$, is described by the equation (see \cite{ABN11})
\begin{align}\label{eqn:diode}
	\begin{cases}
		-\varepsilon^{2}\psi_{E}^{\prime\prime}(x) + V(x) \psi_{E}(x) = E\psi_{E}(x) \Comma \quad x\in(0,1) \Comma \\
		\psi_{E}^{\prime}(0) + \ii k(0)\psi_{E}(0) = 0 \Comma \\
		\psi_{E}^{\prime}(1) - \ii k(1)\psi_{E}(1) = -2\ii k(1) \Comma
	\end{cases}
\end{align}
where $k(x):=\varepsilon^{-1}\sqrt{E-V(x)}$ is the wave vector and $V$ denotes the electrostatic potential. Note that our assumption $a(x) \geq \tau > 0$ implies $E>V(x)$, so the solution $\psi_{E}$ becomes oscillatory. Then, one is often interested in macroscopic quantities like the electron density $n$ and the current density $j$, which are given by
\begin{align} \label{eq:macroscopic}
	n(x)=\int_{0}^{\infty}f(k)\vert \psi_{E(k)}(x)\vert^{2} \,\mathrm{d}k \Comma \quad j(x)=\varepsilon\int_{0}^{\infty}f(k)\Im\left(\overline{\psi_{E(k)}(x)}\psi^{\prime}_{E(k)}(x)\right) \,\mathrm{d}k \Comma
\end{align}
where $f$ represents the injection statistics of the electrons. Here, $\Im(\cdot)$ denotes the imaginary part and $E(k)=\varepsilon^{2}k^{2}+V$ means the energy for a given wave vector $k$. In order to compute the quantities (\ref{eq:macroscopic}), the Schrödinger equation (\ref{eqn:diode}) has to be solved many times, as a fine grid in $E(k)$ is needed. Hence, efficient methods for the solution of (\ref{eqn:diode}) are of great interest in such applications. Instead of solving the boundary value problem (\ref{eqn:diode}) directly, one can also solve the equivalent initial value problem, which results if equation (\ref{schroedinger}) on the interval $(0,1)$ is augmented with the initial values $\varphi(0)=\varphi_{0}=1$ and $\varepsilon\varphi^{\prime}(0)=\varphi_{0}^{\prime}=-\ii \sqrt{a(0)}$ with $a(x)=E-V(x)$. The solution $\varphi$ of this initial value problem and that one of problem (\ref{eqn:diode}) are then related by
\begin{align}
	\psi_{E}(x)=-\frac{2\ii k(1)}{\varphi^{\prime}(1)-\ii k(1)\varphi(1)}\varphi(x)\Comma \nonumber
\end{align}
according to \cite{ABN11}. Indeed, the method proposed in this paper will deal only with initial value problems for the Schrödinger equation (\ref{schroedinger}), but through this equivalence it is equally suitable for solving problem (\ref{eqn:diode}). %\textcolor{red}{Since in \cite{ABN11} an equivalent initial value problem to the boundary value problem (\ref{eqn:diode}) was introduced, and the numerical scheme in Section \ref{chap:WKBMM} actually deals with an initial value problem for the Schrödinger equation (\ref{schroedinger}), the method proposed in this paper is suitable for solving problem (\ref{eqn:diode}).}

In \cite{ABN11,JL03,LJL05}, efficient and accurate WKB-based (named after the physicists Wentzel, Kramers, Brillouin; cf.\ \cite{LL85}) numerical schemes have been developed for (\ref{schroedinger}) in the oscillatory regime. By transforming out the dominant oscillations, they allow to compute a solution using a coarse spatial grid with step size $h>\lambda$. In fact, the grid limitation can there be reduced to at least $h=\mathcal{O}(\sqrt{\varepsilon})$. Now, this work adds on top of the algorithm from \cite{ABN11} an adaptive step size control as well as a switching mechanism. This allows the algorithm to switch to a standard ODE method (e.g., Runge-Kutta) during the computation in order to avoid technical or efficiency problems in regions where the coefficient function $a(x)$ is very small or indeed equal to zero. We recall that the WKB-approximation is not valid close to turning points, i.e.\ where $a(x)=0$. A switching mechanism was also used in \cite{HLH16}, where the authors presented another WKB-based numerical scheme for the initial value problem corresponding to (\ref{schroedinger}). Therefore, one goal of this work is to compare numerical results from our method with the results given by the method from \cite{HLH16}, by considering two examples where the analytical solution is known.

Since numerical methods for oscillatory problems is an active field of research, let us mention some references that are intended for more general oscillatory problems, and hence also include adequate methods for (\ref{schroedinger}): first the two monographs \cite{WW18,HLW06}. The \textit{adiabatic integrators} of \cite[\S XIV]{HLW06} are in fact closely related to a zeroth order WKB-approximation, see (\ref{wkb_ansatz})-(\ref{wkb_basis0}), below. Concerning (highly oscillatory) \textit{Hamiltonian boundary value methods} we cite \cite{BIMR19,ABI20}.

This paper is organized as follows: In Section \ref{chap:WKBMM} we give a short review of the second order (w.r.t.\ $\varepsilon$) WKB-marching method from \cite{ABN11}. Section \ref{chap:stepsize} then describes the adaptive step size control algorithm as well as the switching mechanism. In Section \ref{chap:RKWKB} we recap the Runge-Kutta-WKB method from \cite{HLH16} and point out the difference between their step size control and the one used in this paper. In Section \ref{chap:numerical_results} we present numerical investigations on the error estimators of our algorithm as well as a comparison of the numerical results of our method and the method from \cite{HLH16}. We then conclude in Section \ref{chap:conclusion}.
\section{The WKB-marching method}\label{chap:WKBMM}
We aim at solving the Schrödinger equation (\ref{schroedinger}), augmented with the initial conditions
\begin{align} \label{initial_condition}
	\varphi(x_{0})=\varphi_{0}\Comma \quad \varepsilon\varphi^{\prime}(x_{0})=\varphi_{0}^{\prime}
\end{align}
with some $x_{0}\in\mathbb{R}$. First we shall review the basics of the second order (w.r.t.\ $\varepsilon$) WKB-marching method from \cite{ABN11} with focus on the algorithm. The motivation for this method was the construction of a numerical scheme that is uniformly correct in $\varepsilon$ and sometimes even asymptotically correct, i.e.\ the numerical error goes to zero with $\varepsilon\to 0$ while the grid size $h$ remains constant. For further details see \cite{ABN11}. The method consists of two parts:
\begin{enumerate}
	\item Analytic pre-processing of (\ref{schroedinger}) by transforming the equation into a smoother (i.e.\ less oscillatory) problem that can be solved accurately and efficiently on a coarse grid.
	\item Obtaining a numerical scheme by discretization of the smoother problem.
\end{enumerate}
\textbf{Analytic pre-processing.}
The well-known WKB-approximation (cf.\ \cite{LL85}) for the oscillatory regime where $a(x)\geq \tau$ for some $\tau > 0$, is based on inserting the ansatz
\begin{align} \label{wkb_ansatz}
	\varphi(x)\sim \exp\left(\frac{1}{\varepsilon}\sum_{p=0}^{\infty}\varepsilon^{p}\phi_{p}(x)\right)
\end{align}
into equation (\ref{schroedinger}). After a comparison of $\varepsilon$-powers one obtains the first three functions $\phi_{p}(x)$ as
\begin{align}
	\phi_{0}(x)&= \pm \ii \int_{x_{0}}^{x} \sqrt{a(y)} \, \mathrm{d}y \Comma \label{wkb_basis0}\\
	\phi_{1}(x)&= \ln(a(x)^{-\frac{1}{4}}) \Comma \label{wkb_basis1}\\
	\phi_{2}(x)&= \mp \ii \int_{x_{0}}^{x} b(y) \, \mathrm{d}y \Comma \quad b(x):= -\frac{1}{2a(x)^{\frac{1}{4}}}\left(a(x)^{-\frac{1}{4}}\right)^{\prime\prime} \period \label{wkb_basis2}
\end{align}
Here the symbols $\pm$ and $\mp$ in (\ref{wkb_basis0}) and (\ref{wkb_basis2}) correspond to the fact that there is always a pair of approximate solutions to the Schrödinger equation (\ref{schroedinger}), by analogy to the two fundamental solutions of (\ref{schroedinger}). Hence the general solution is then a linear combination of the two. Therefore, a truncation of the sum (\ref{wkb_ansatz}) after $p = 2$ leads to the second order (w.r.t.\ $\varepsilon$) asymptotic WKB-approximation of $\varphi(x)$:
\begin{align} \label{2ndorderwkb}
	\varphi_{2}(x) = c_{1}\frac{\exp\left(\frac{\ii}{\varepsilon} \phi^{\varepsilon}(x)\right)}{a(x)^{\frac{1}{4}}} + c_{2}\frac{\exp\left(-\frac{\ii}{\varepsilon}\phi^{\varepsilon}(x)\right)}{a(x)^{\frac{1}{4}}} \Comma
\end{align}
with constants $c_{1},c_{2}\in\mathbb{C}$ to be determined from initial or boundary conditions, and the phase is
\begin{align} \label{phase}
	\phi^{\varepsilon}(x):= \int_{x_{0}}^{x}\left(\sqrt{a(y)}-\varepsilon^{2}b(y)\right) \, \mathrm{d}y \period
\end{align}
In the WKB-marching method of \cite{ABN11} this second order WKB\hyp{}approximation is used to transform (\ref{schroedinger}) into a smoother problem. To clarify our terminology, we point out that this method (as well as the Runge-Kutta-WKB method of Section \ref{chap:RKWKB}) has both a \textit{WKB-order} (w.r.t.\ $\varepsilon$; referring to the used cut-off in the asymptotic expansion (\ref{wkb_ansatz})) and a \textit{numerical order} (w.r.t.\ the step size $h$; referring to the convergence order). Firstly, using the notation
\begin{align} \label{Utrafo}
	U(x)=\begin{pmatrix}u_{1}(x)\\u_{2}(x)\end{pmatrix}:=\begin{pmatrix}a(x)^{\frac{1}{4}}\varphi(x)\\\frac{\varepsilon\left(a(x)^{\frac{1}{4}}\varphi(x)\right)^{\prime}}{\sqrt{a(x)}}\end{pmatrix} \Comma
\end{align}
the second order differential equation (\ref{schroedinger}) with the initial conditions (\ref{initial_condition}) can be reformulated as a system of first order differential equations:
\begin{align} \label{Usystem}
	\begin{cases}
		U^{\prime}(x) = \left[\frac{1}{\varepsilon}\mathbf{A}_{0}(x)+\varepsilon \mathbf{A}_{1}(x)\right]U(x) \Comma \quad x>x_{0}\Comma\\
		U(x_{0})=U_{I} \period
	\end{cases}
\end{align}
Here, the two matrices $\mathbf{A}_{0}$ and $\mathbf{A}_{1}$ are given by
\begin{align}
	\mathbf{A}_{0}(x):=\sqrt{a(x)}\begin{pmatrix}0 & 1 \\ -1 & 0\end{pmatrix}\semicolon\quad \mathbf{A}_{1}(x):=\begin{pmatrix}0 & 0 \\ 2b(x) & 0\end{pmatrix} \nonumber \period
\end{align}
Then, the first order system (\ref{Usystem}) for $U(x)$ is transformed by the change of variables
\begin{align} \label{Z-trafo}
	Z(x)=\begin{pmatrix}z_{1}(x) \\ z_{2}(x)\end{pmatrix}:=\exp\left(-\frac{\ii}{\varepsilon}\mathbf{\Phi}^{\varepsilon}(x)\right)\mathbf{P}U(x) \Comma
\end{align}
with the two matrices
\begin{align}
	\mathbf{P}:=\frac{1}{\sqrt{2}}\begin{pmatrix}\ii & 1 \\ 1 & \ii\end{pmatrix}\semicolon \quad \mathbf{\Phi}^{\varepsilon}(x):=\begin{pmatrix}\phi^{\varepsilon}(x) & 0 \\ 0 & -\phi^{\varepsilon}(x) \end{pmatrix} \nonumber \Comma
\end{align}
where $\phi^{\varepsilon}$ is the phase function defined in (\ref{phase}).
%leads to
%\begin{align}
%\begin{cases}
%Y^{\prime}(x) = \frac{\ii}{\varepsilon}D^{\varepsilon}(x)Y(x)+\varepsilon N(x)Y(x) \Comma \quad x\in\\
%Y(x_{0})=Y_{I}
%\end{cases}
%\end{align}
%with
%\begin{align}
%D^{\varepsilon}(x)=\begin{pmatrix}\sqrt{a(x)}-\varepsilon^{2}\beta(x) & 0 \\ 0 & -\sqrt{a(x)}+\varepsilon^{2}\beta(x)\end{pmatrix}\semicolon\quad N(x):=\begin{pmatrix}0 & \beta(x) \\ \beta(x) & 0\end{pmatrix}\period
%\end{align}
%A final transformation for the elimination of leading oscillations by defining the diagonal matrix
%\begin{align}
%	\Phi^{\varepsilon}(x):=\begin{pmatrix}\phi^{\varepsilon}(x) & 0 \\ 0 & -\phi^{\varepsilon}(x) \end{pmatrix}
%\end{align}
%and a change of unknown
%\begin{align}
%	Z(x)=\begin{pmatrix}z_{1}(x) \\ z_{2}(x)\end{pmatrix}:=\exp\left(-\frac{\ii}{\varepsilon}\Phi^{\varepsilon}(x)\right)Y(x) \comma
%\end{align}
This leads to the system
\begin{align} \label{Zsystem}
	\begin{cases}
		Z^{\prime}(x) = \varepsilon \mathbf{N}^{\varepsilon}(x)Z(x) \Comma \quad x>x_{0}\Comma\\
		Z(x_{0})=Z_{I}=\mathbf{P}U_{I} \Comma
	\end{cases}
\end{align}
where $\mathbf{N}^{\varepsilon}(x)$ is a (Hermitian) matrix with only off-diagonal non-zero entries:
\begin{align}
	N_{1,2}^{\varepsilon}(x)=b(x)\operatorname{e}^{-\frac{2\ii}{\varepsilon}\phi^{\varepsilon}(x)}\Comma\quad N_{2,1}^{\varepsilon}(x)=b(x)\operatorname{e}^{\frac{2\ii}{\varepsilon}\phi^{\varepsilon}(x)} \period \nonumber
\end{align}
%
%\begin{align}
%	N^{\varepsilon}(x):=\exp\left(-\frac{\ii}{\varepsilon}\Phi^{\varepsilon}(x)\right)N(x)\exp\left(\frac{\ii}{\varepsilon}\Phi^{\varepsilon}(x)\right) \period
%\end{align}
Since the transformation (\ref{Z-trafo}) eliminated the dominant oscillations, the system (\ref{Zsystem}) can be solved numerically on a coarse grid $\{x_{n}, n\in\mathbb{N}_{0}\}$. Then the original solution can be recovered by the inverse transformation
\begin{align}
	U(x)=\mathbf{P}^{-1}\exp\left(\frac{\ii}{\varepsilon}\mathbf{\Phi}^{\varepsilon}(x)\right)Z(x) \period \nonumber
\end{align}
It should be noted that, throughout the whole transformation from $U(x)$ to $Z(x)$, the phase integral $\phi^{\varepsilon}(x)$ is assumed to be known exactly. For a generalization using a spectral method to numerically compute the phase see \cite{AKU19}.
\medskip

\noindent \textbf{Numerical scheme.} The derivation of the second order (in $h$) scheme for (\ref{Zsystem}) is obtained via the second order Picard approximation
\begin{align} \label{picard}
	\widetilde{Z}_{n+1}:=\widetilde{Z}_{n}+\varepsilon\int_{x_{n}}^{x_{n+1}}\mathbf{N}^{\varepsilon}(x) \, \mathrm{d}x \widetilde{Z}_{n} + \varepsilon^{2}\int_{x_{n}}^{x_{n+1}}\mathbf{N}^{\varepsilon}(x)\int_{x_{n}}^{x}\mathbf{N}^{\varepsilon}(y) \, \mathrm{d}y\mathrm{d}x \widetilde{Z}_{n}\period
\end{align}
Since the entries of $\mathbf{N}^{\varepsilon}(x)$ are highly oscillatory, (\ref{picard}) involves (iterated) oscillatory integrals. With $\phi^{\varepsilon}$ assumed to be known exactly, they are then approximated using similar techniques as the \textit{asymptotic method} in \cite{INO06}. The first order (in $h$) scheme for (\ref{Zsystem}) is derived by only taking into account the first two terms from (\ref{picard}). For both of these schemes we introduce the following notations:
\begin{align}
	&b_{0}(y):=\frac{b(y)}{2\left(\sqrt{a(y)}-\varepsilon^{2}b(y)\right)}\semicolon\quad b_{k+1}(y):=\frac{1}{2\left(\phi^{\varepsilon}(y)\right)^{\prime}}\frac{\mathrm{d}b_{k}}{\mathrm{d}y}(y)\Comma \quad k=0,1,2\semicolon\nonumber\\
	&h_{1}(y):=\operatorname{e}^{\ii y}-1\semicolon\quad h_{2}(y):=\operatorname{e}^{\ii y}-1-\ii y\period \nonumber
\end{align}
Further, let $\{x_{0},x_{1},...,x_{N}\}$ be a grid we want to compute the solution on, and $h:=\max_{1\leq n\leq N}\vert x_{n}-x_{n-1}\vert$ be the step size. Then both schemes read as follows:
\medskip

\noindent \underline{First order scheme:}
Let $Z_{0}:=Z_{I}$ be the initial condition and let $n=0,...,N-1$. Then the algorithm updates as
\begin{align}
	Z_{n+1}=\left(I+\mathbf{A}_{n}^{1}\right)Z_{n} \Comma \label{scheme:first_order}
\end{align}
with the (Hermitian) matrix
\begin{align}
	&\mathbf{A}_{n}^{1}:=\varepsilon^{3}b_{1}(x_{n+1})\begin{pmatrix}0 & \operatorname{e}^{-\frac{2\ii}{\varepsilon}\phi^{\varepsilon}(x_{n})}h_{1}\left(-\frac{2}{\varepsilon}s_{n}\right) \\ \operatorname{e}^{\frac{2\ii}{\varepsilon}\phi^{\varepsilon}(x_{n})}h_{1}\left(\frac{2}{\varepsilon}s_{n}\right) & 0\end{pmatrix} \nonumber \\
	&-\ii \varepsilon^{2}\begin{pmatrix}0 & b_{0}(x_{n})\operatorname{e}^{-\frac{2\ii}{\varepsilon}\phi^{\varepsilon}(x_{n})}-b_{0}(x_{n+1})\operatorname{e}^{-\frac{2\ii}{\varepsilon}\phi^{\varepsilon}(x_{n+1})}\\ b_{0}(x_{n+1})\operatorname{e}^{\frac{2\ii}{\varepsilon}\phi^{\varepsilon}(x_{n+1})}-b_{0}(x_{n})\operatorname{e}^{\frac{2\ii}{\varepsilon}\phi^{\varepsilon}(x_{n})} & 0\end{pmatrix}\Comma \nonumber
\end{align}
and the phase increments
\begin{align}
	s_{n}:=\phi^{\varepsilon}(x_{n+1})-\phi^{\varepsilon}(x_{n})\period \nonumber
\end{align}
\medskip

\noindent \underline{Second order scheme:} Let $Z_{0}:=Z_{I}$ be the initial condition and let $n=0,...,N-1$. Then the algorithm updates as
\begin{align}
	Z_{n+1}=\left(I+\mathbf{A}_{mod,n}^{1}+\mathbf{A}_{n}^{2}\right)Z_{n} \Comma \label{scheme:second_order}
\end{align}
with the (Hermitian) matrix
\begin{align}
	&\mathbf{A}_{mod,n}^{1}:=\nonumber\\
	&-\ii\varepsilon^{2}\begin{pmatrix}0 & b_{0}(x_{n})\operatorname{e}^{-\frac{2\ii}{\varepsilon}\phi^{\varepsilon}(x_{n})}-b_{0}(x_{n+1})\operatorname{e}^{-\frac{2\ii}{\varepsilon}\phi^{\varepsilon}(x_{n+1})} \\ b_{0}(x_{n+1})\operatorname{e}^{\frac{2\ii}{\varepsilon}\phi^{\varepsilon}(x_{n+1})}-b_{0}(x_{n})\operatorname{e}^{\frac{2\ii}{\varepsilon}\phi^{\varepsilon}(x_{n})} & 0 \end{pmatrix} \nonumber \\
	&+\varepsilon^{3} \begin{pmatrix} 0 & b_{1}(x_{n+1})\operatorname{e}^{-\frac{2\ii}{\varepsilon}\phi^{\varepsilon}(x_{n+1})}-b_{1}(x_{n})\operatorname{e}^{-\frac{2\ii}{\varepsilon}\phi^{\varepsilon}(x_{n})} \\ b_{1}(x_{n+1})\operatorname{e}^{\frac{2\ii}{\varepsilon}\phi^{\varepsilon}(x_{n+1})}-b_{1}(x_{n})\operatorname{e}^{\frac{2\ii}{\varepsilon}\phi^{\varepsilon}(x_{n})} & 0\end{pmatrix} \nonumber \\
	&+\ii\varepsilon^{4}b_{2}(x_{n+1})\begin{pmatrix} 0 & -\operatorname{e}^{-\frac{2\ii}{\varepsilon}\phi^{\varepsilon}(x_{n})}h_{1}\left(-\frac{2}{\varepsilon}s_{n}\right) \\ \operatorname{e}^{\frac{2\ii}{\varepsilon}\phi^{\varepsilon}(x_{n})}h_{1}\left(\frac{2}{\varepsilon}s_{n}\right) & 0 \end{pmatrix} \nonumber\\ &-\varepsilon^{5}b_{3}(x_{n+1}) \begin{pmatrix} 0 & \operatorname{e}^{-\frac{2\ii}{\varepsilon}\phi^{\varepsilon}(x_{n})}h_{2}\left(-\frac{2}{\varepsilon}s_{n}\right) \\ \operatorname{e}^{\frac{2\ii}{\varepsilon}\phi^{\varepsilon}(x_{n})}h_{2}\left(\frac{2}{\varepsilon}s_{n}\right) & 0\end{pmatrix}\Comma \nonumber
\end{align}
and the diagonal matrix
\begin{align}
	\mathbf{A}_{n}^{2}:= &-\ii\varepsilon^{3}\left(x_{n+1}-x_{n}\right)\frac{b(x_{n+1})b_{0}(x_{n+1})+b(x_{n})b_{0}(x_{n})}{2}\begin{pmatrix}1 & 0 \\ 0 & -1\end{pmatrix} \nonumber\\
	&-\varepsilon^{4}b_{0}(x_{n})b_{0}(x_{n+1})\begin{pmatrix}h_{1}\left(-\frac{2}{\varepsilon}s_{n}\right) & 0 \\
		0 & h_{1}\left(\frac{2}{\varepsilon}s_{n}\right) \end{pmatrix}\nonumber\\
	&+\varepsilon^{5}b_{1}(x_{n+1})\left[b_{0}(x_{n})-b_{0}(x_{n+1})\right]\begin{pmatrix}h_{2}\left(-\frac{2}{\varepsilon}s_{n}\right) & 0 \\
		0 & -h_{2}\left(\frac{2}{\varepsilon}s_{n}\right)\end{pmatrix}\period \nonumber
\end{align}
\section{Step size control and switching mechanism}
%\section{\textcolor{red}{Step size control and coupling mechanism}}
\label{chap:stepsize}
The WKB scheme is efficient in the highly oscillatory regime, but not applicable close to turning points, i.e.\ zeros of $a(x)$, see \cite{AD20}. This is evident already from the transformation (\ref{Utrafo}), which does not make sense when $a(x)\leq 0$. For mixed problems, e.g., the Airy equation on $\mathbb{R}_{0}^{+}$ (see Section \ref{airy_section}), which has a turning point at $x=0$, it is therefore convenient to couple  two different methods: a method for highly oscillatory ODEs (e.g., WKB-based) away from the turning point, and a standard ODE method (e.g., Runge-Kutta) close to the turning point, where the solution is smooth anyhow. Here, we choose the well-known Runge-Kutta-Fehlberg 4(5) (RKF45) scheme (cf.\ \cite{HNW00}) as the standard ODE method. The latter method will be applied directly on equation (\ref{schroedinger}) and not on (\ref{Zsystem}), since the WKB-transformation (\ref{Utrafo})-(\ref{Zsystem}) is not permitted at turning points. The exact switching mechanism as well as the introduction of an adaptive step size control to the algorithm will be described in the two following subsections.
\subsection{The adaptive step size controller}
In order to compute the solution efficiently, an adaptive step size controller, based on an estimator for the local truncation error, will be added to the numerical methods. This control allows the step size to increase or decrease while aiming to keep the error estimator as close as possible within a given error tolerance. %Furthermore, in each step, two solutions are computed, namely, the first one  based on the WKB-marching method and the second one via a Runge-Kutta-Fehlberg 4(5) scheme. The step size control then not only chooses the respective method with the step size expected to be larger in the next step, but also aims to ensure that the numerical error, which is done by going one step forward, to be within a given relative error tolerance $\operatorname{RTol}$.
%The Switching mechanism is used for efficiency reasons. The RKF 4(5) scheme is more efficient in the low oscillating regime whereas the WKB-marching method becomes even more efficient (in terms of step size) the higher the oscillations are.
%The step size control is using an estimator for the local truncation error at each step.
To illustrate the functionality of this step size controller, we shall consider a numerical scheme of order $k$. We are then able to apply this step size control individually to the different methods mentioned above.\\
Let $Y_{n}^{(k)}=(\varphi_{n}^{(k)}, \varphi_{n}^{\prime(k)})$ and $Y_{n}^{(k+1)}=(\varphi_{n}^{(k+1)}, \varphi_{n}^{\prime(k+1)})$ be two numerical solutions of order $k$ and $k+1$ to approximate the exact solution $Y(x_{n})=(\varphi(x_{n}), \varphi^{\prime}(x_{n}))$ of the initial value problem (\ref{schroedinger}), (\ref{initial_condition}). E.g., one could choose the WKB schemes (\ref{scheme:first_order}) and (\ref{scheme:second_order}) of $h$-order 1 and 2, respectively. Next we want to decide whether to accept the numerical solution at $x_{n}$ (typically the more accurate solution $Y_{n}^{(k+1)}$) or rather to retry the computation with a modified step size. To this end we define the estimator for the local truncation error as
\begin{align} \label{wkb_est}
\operatorname{est}_{n}:=\lVert Y_{n}^{(k)}-Y_{n}^{(k+1)}\rVert_{\infty} \period
\end{align}
%
%At this point it is not clear if this is an adequate choice for the estimator, which is why we will compare its values with the real local truncation errors using reference solutions. Also it would make sense to use an error estimator based on the not retransformed $Z$-variable, however, for the sake of comparability we use the definition above since there is no $Z$-variable in the RKWKB method. \textcolor{blue}{Zur Überlegung bzgl. einer kurzen Anmerkung, dass man den Schätzer auch auf $Z$-Level definieren könnte: Die Transformation von $U$ nach $Z$ ist zwar unitär, jedoch ist $U$ nicht einfach nur ein skaliertes $(\varphi, \varphi^{\prime})$, denn es ist $u_{2}=\varepsilon\left(\frac{1}{4}a(x)^{-\frac{5}{4}}a^{\prime}(x)\varphi(x)+a(x)^{-\frac{1}{4}}\varphi^{\prime}(x)\right)$. Beim Airy-Beispiel hatte man damals keinen signifikanten Unterschied beim Verwenden eines Schätzers auf $Z$-Level gesehen, vermutlich weil der erste Summand von $u_{2}$ für großes $x$ immer unwichiger wird. Für das allgemeine Problem (\ref{schroedinger}) ist das aber nicht notwendigerweise so.}
%Note also that $U_{n}$ is connected to $Z_{n}$ just by an unitary transformation and therefore the error estimator could easily be defined on the $Z$-level of the solution, since $y_{WKB}^{(.)}$ is essentially a scaled version of $U$\\
Let $h_{n,trial}:=x_{n}-x_{n-1}$ be the (trial) step size which was used to compute the solutions at the current step $n$. We then use the common approach of varying the step size via the multiplicative control
\begin{align}
	h_{new}&:=\theta_{n} \cdot h_{n,trial} \period \nonumber
\end{align}
Here, we choose the factor $\theta_{n}$ to be based on the so-called \textit{elementary controller} (e.g., see \cite{Soed04}). Additionally, we introduce limitations in such way that the step size controller responds ``smoothly'' to abrupt changes in the solution behaviour, that is, the ratio between two consecutive step sizes should not be exorbitantly large or small. That said, we choose the factor similar to \cite[p. 310]{But08} as
\begin{align}
\theta_{n}:=\max\left(0.5,\min\left(2,0.9\left(\frac{\operatorname{ATol}+\operatorname{RTol}\cdot\lVert Y_{n}^{(k+1)}\rVert_{\infty}}{\operatorname{est}_{n}}\right)^{\frac{1}{k+1}}\right)\right)\Comma \label{theta}
\end{align}
where $\operatorname{ATol}=\eta\cdot\operatorname{Tol}$ and $\operatorname{RTol}=\operatorname{Tol}$ are absolute and relative error tolerances for a given master tolerance $\operatorname{Tol}$, the values $0.5$ and $2$ are design parameters that limit the ratio of two consecutive step sizes from below and above, and $0.9$ is a common safety factor for increasing the probability of the next step to be accepted. Here, $\eta$ is a scaling factor representing the gradual switch-over between absolute and relative errors, depending on the behaviour of the solution. This is because for $\lVert Y_{n}^{(k+1)}\rVert_{\infty}\to 0$ the $\operatorname{ATol}$ term in the numerator in (\ref{theta}) is dominating, whereas for large values of $\lVert Y_{n}^{(k+1)}\rVert_{\infty}$ the $\operatorname{RTol}$ term is dominating. For our numerical simulations in Section \ref{chap:numerical_results} we choose $\eta=10^{-2}$. If $\operatorname{ATol}+\operatorname{RTol}\cdot\lVert Y_{n}^{(k+1)}\rVert_{\infty}\geq\operatorname{est}_{n}$, we accept the $n$-th step with the step size defined as $h_{n}:=h_{n,trial}$ and define the trial step size for the next step as
\begin{align}
	h_{n+1,trial}&:=h_{new} \period \nonumber
\end{align}
However, if $\operatorname{ATol}+\operatorname{RTol}\cdot\lVert Y_{n}^{(k+1)}\rVert_{\infty}<\operatorname{est}_{n}$, the $n$-th step gets rejected and a reattempt is done with the smaller (trial) step size $h_{n,trial}$ by updating its value as
\begin{align}
	h_{n,trial}\to h_{new} \period \nonumber
\end{align}
Since such an acceptance criterion is based on aiming to maintain the local error in each step as close as possible to the given error tolerances it is often referred to as \textit{error per step (EPS)} control. In practice, using EPS control one can hope to achieve a global truncation error proportional to $\operatorname{Tol}^{k/(k+1)}$ (e.g., see \cite[p. 311]{But08}).
\subsection{The switching mechanism}
As already mentioned above, the algorithm shall automatically switch between two numerical methods, the WKB method in the oscillatory regime and another method, which is valid close to turning points. To realize this dynamical switching mechanism we follow a similar strategy as in \cite{AHLH20}. To illustrate this procedure we now consider two numerical schemes of order $k^{(1)}$ and $k^{(2)}$, where the superscripts (1) and (2) correspond to the two schemes. In each step, the adaptive step size algorithm from the previous section is applied to both schemes individually up to the definition of the quantities $\theta_{n}^{(1)}$ and $\theta_{n}^{(2)}$, i.e., we just evaluate (\ref{wkb_est}) and (\ref{theta}). Then, after checking the acceptance criterion $\operatorname{ATol}+\operatorname{RTol}\cdot\lVert Y_{n}^{(k^{(i)}+1)}\rVert_{\infty}\geq\operatorname{est}_{n}^{(i)}$ for each scheme $(i)\in\{(1),(2)\}$, the switching mechanism intervenes and it selects the acceptable numerical method that yields the larger value of $\theta_{n}^{(i)}$, for $i=1,2$, hence yielding the larger proposed new step size. We thus favor the method with the smaller error estimator, discounted by its respective order $k^{(i)}$. More precisely, we define
\begin{align}
\Theta_{n}:=
	\begin{cases}
	\theta_{n}^{(1)} \Comma &(1)\text{ accepted, }(2)\text{ rejected}\\
	\theta_{n}^{(2)} \Comma &(1)\text{ rejected, }(2)\text{ accepted}\\
	\max\left(\theta_{n}^{(1)},\theta_{n}^{(2)}\right) \Comma &\text{otherwise}
	\end{cases} \nonumber
\end{align}
and store the information on the method of choice in that $n$-th step. Through this procedure the algorithm does not only use the error estimator to find the next step size, but also to decide between the two methods.

If at least one method was accepted, the algorithm sets
\begin{align}
	h_{n}&:=h_{n,trial} \Comma \nonumber \\
	h_{n+1,trial}&:=\Theta_{n}\cdot h_{n,trial} \period \nonumber
\end{align}
Otherwise a reattempt is done with the smaller (trial) step size $h_{n,trial}$ by updating its value as
\begin{align}
	h_{n,trial}\to \Theta_{n} \cdot h_{n,trial} \period \nonumber
\end{align}

We remark that the coupling of two methods, as presented above, could incur additional computational costs, since in every step both methods have to be applied in order to compute $\Theta_{n}$. However, in our case the WKB method (for highly oscillatory regions) and the standard ODE solver (for smoother regions) complement each other very well, yielding better results concerning accuracy as well as overall efficiency (see also Figures \ref{plot:Airy_RKWKBMM_ratios}-\ref{plot:Airy_RKWKB_ratios}, \ref{plot:Airy_WKBMM_RKWKBmod_RKWKB_RKF45_CPUtime_vs_err} and \ref{plot:PCF_WKBMM_RKWKBmod_RKWKB_RKF45_CPUtime_vs_err}).
\section{The Runge-Kutta-WKB method}\label{chap:RKWKB}
In this section we give a short review of the Runge-Kutta-WKB (RKWKB) method presented in \cite{HLH16} for the initial value problem (\ref{schroedinger}), (\ref{initial_condition}). The method is based on a dynamic switching mechanism between a standard Runge-Kutta scheme and a new stepping procedure that uses the WKB-ansatz (\ref{wkb_ansatz}) as an approximation of the true solution. This stepping procedure reads as follows:
\begin{align}
	\varphi_{n+1}&:=\gamma_{+}f_{+}(x_{n}+h_{n})+\gamma_{-}f_{-}(x_{n}+h_{n}) \Comma \label{stepping_start}\\
	\varphi^{\prime}_{n+1}&:=\delta_{+}f^{\prime}_{+}(x_{n}+h_{n})+\delta_{-}f^{\prime}_{-}(x_{n}+h_{n}) \Comma \label{stepping_start2}\\
	x_{n+1}&:=x_{n}+h_{n} \label{stepping_start3}\Comma
\end{align}
where
\begin{align}
	\gamma_{\pm}&:=\frac{\varphi^{\prime}_{n}f_{\mp}(x_{n})-\varphi_{n}f^{\prime}_{\mp}(x_{n})}{f^{\prime}_{\pm}(x_{n})f_{\mp}(x_{n})-f^{\prime}_{\mp}(x_{n})f_{\pm}(x_{n})} \Comma \\
	\delta_{\pm}&:=\frac{\varphi^{\prime\prime}_{n}f^{\prime}_{\mp}(x_{n})-\varphi^{\prime}_{n}f^{\prime\prime}_{\mp}(x_{n})}{f^{\prime\prime}_{\pm}(x_{n})f^{\prime}_{\mp}(x_{n})-f^{\prime\prime}_{\mp}(x_{n})f^{\prime}_{\pm}(x_{n})} \Comma \label{stepping_end}
\end{align}
and $\varphi_{n}^{\prime\prime}$ is computed from $\varphi_{n}$ and equation (\ref{schroedinger}).
Here, $f_{\pm}$ are chosen as WKB-approximations (\ref{wkb_ansatz}) of some finite order. For instance, one gets the second (WKB-)order method by setting
\begin{align}
	f_{\pm}(x)&:=\frac{\exp\left(\pm\frac{\ii}{\varepsilon}\phi^{\varepsilon}(x)\right)}{a(x)^{\frac{1}{4}}} \period \nonumber
\end{align}
Note that this choice is equal to (\ref{2ndorderwkb}), but one can easily choose $f_{\pm}$ with higher (WKB-)orders. However, the stepping procedure (\ref{stepping_start})-(\ref{stepping_end}) is always a first order (in $h$) numerical method. This holds because the coefficients $\gamma_{\pm}$ and $\delta_{\pm}$ are chosen in such way that one finds
\begin{align}
	\varphi_{n+1} &= \varphi_{n} + \varphi_{n}^{\prime}h + \mathcal{O}(h^{2}) \Comma \nonumber \\
	\varphi_{n+1}^{\prime} &= \varphi_{n}^{\prime} + \varphi_{n}^{\prime\prime}h + \mathcal{O}(h^{2}) \Comma \nonumber
\end{align}
from (\ref{stepping_start})-(\ref{stepping_start3}), as stated in \cite{HLH16}.

It is also worth noting that in \cite{HLH16} the authors use a slightly different dynamic switching mechanism and a different step size control compared to the algorithm presented in Section \ref{chap:stepsize}. Firstly, their estimator of the relative error within the WKB-stepping procedure uses the difference between two numerical solutions $\varphi_{n}^{(k)}$ and $\varphi_{n}^{(k+1)}$ of different WKB-orders instead of different $h$-orders, simply since they do not have schemes of two different $h$-orders at their disposal. The algorithm then decides between a WKB step and a RK step by choosing the method with the smaller error estimator. In their more recent paper \cite{AHLH20} the authors use a similar step size control and switching mechanism as presented in Section \ref{chap:stepsize}, but they do not limit the ratio of two consecutive step sizes and provide no option for controlling both absolute and relative errors as done in (\ref{theta}).

\medskip
%\textcolor{red}{It should be noted that within this method the WKB-order does not necessarily equal the order w.r.t. $h$. Error estimator based on two numerical solutions of different WKB order is just based on two solutions of same $h$-order. Whereas WKB-marching method: $p_{WKB}=2$ implicates second order scheme in $h$!}\\
The goal of the following section is to compare numerical results from the WKB-marching method to results one gets using the RKWKB method. To both methods we will apply (exactly) the step size control and switching mechanism from Section \ref{chap:stepsize}, for the sake of comparability. Since the WKB-stepping procedure of the RKWKB method is always of first order w.r.t.\ $h$, the definition of the error estimator (\ref{wkb_est}) does not make sense. Therefore we shall use two different WKB orders (instead of $h$-orders) to be able to compute the error estimator in this case, as also done in \cite{HLH16}.
Further, since our algorithm consists of a different step size update formula and switching criterion, we will call this modified RKWKB method simply RKWKBmod. Since this modification may produce slightly different numerical results compared to \cite{HLH16,AHLH20} we will also include the original RKWKB method from \cite{HLH16} into our comparison.
\section{Numerical results: WKB-marching method vs.\ the Runge-Kutta-WKB method}\label{chap:numerical_results}
%
%\textcolor{blue}{Am Ende alle Plots nochmal erstellen und dabei die Referenzlösungen möglichst fein abtasten! Sonst sieht man beim Zoom auf einen Plot vor allem bei hoch oszillierenden Lösungen die Gitterpunkte, was unschön ist.}
In this section we will compare numerical results of the WKB-marching method and of the RKWKB method by applying both algorithms to two examples, where exact analytical solutions are available. The first example corresponds to a linear coefficient function $a(x)$ and is taken from \cite{HLH16,AHLH20}, whereas the second example involves a quadratic function $a(x)$ and appears in \cite{AD20}. In both examples the phase integral (\ref{phase}) in the WKB basis functions (\ref{wkb_basis0})-(\ref{wkb_basis2}) can easily be computed exactly, hence we do not need any numerical integration routine here. By contrast, in \cite{AHLH20} they evaluate the WKB basis functions with a numerical integration routine and therefore get another source of error in their method. %Note also that the computation of the WKB basis functions (\ref{wkb_basis0}) and (\ref{wkb_basis2}) directly corresponds to the computation of the phase (\ref{phase}).
Moreover, we will always use the second order WKB-marching method, since no scheme of higher WKB-order has been developed yet, see \cite{ABN11}.

For clarity of the presentation we shall use in the sequel the following terminology for the methods to be compared:
\begin{table}[H] 
	\centering
	\begin{tabularx}{\textwidth}{|l|X|}
		\hline
		WKB+RKF45 & WKB-marching method (see Section \ref{chap:WKBMM}) + step size and switching algorithms to RKF45 (see Section \ref{chap:stepsize}) \\
		\hline
		RKWKB & original method from \cite{HLH16}: Runge-Kutta-WKB (see Section \ref{chap:RKWKB}) + original step size and switching algorithms (see Section \ref{chap:RKWKB}) \\
		\hline
		RKWKBmod & modified method from \cite{HLH16}: Runge-Kutta-WKB (see Section \ref{chap:RKWKB}) + modified step size and switching algorithms (see Section \ref{chap:stepsize}) \\
		\hline
		RKF45 & Runge-Kutta-Fehlberg 4(5) scheme + step size algorithm (see Section \ref{chap:stepsize}) \\
		\hline
	\end{tabularx}
%	\begin{tabularx}{\textwidth}{ |l|X| }
%		\hline
%		word1 & long definition... \\
%		\hline 
%		word2  & long definition...  \\
%		\hline
%	\end{tabularx}
	\caption{Terminology for the methods to be compared.}
	
	\label{table:methods}
\end{table}
All simulations in this paper are done with \MATLAB version 9.10.0.1669831 (R2021a).
\subsection{First example: Airy function} \label{airy_section}
\begin{figure}[H]
	\centering
	\includegraphics[scale=0.95]{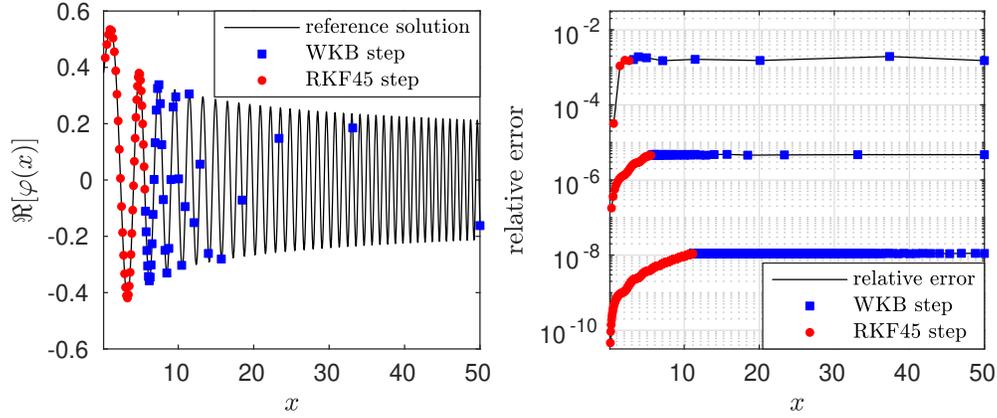}
	\caption{Left: Real part of the numerical solution obtained by using WKB+RKF45 compared to the (exact) reference solution (solid line) for $\operatorname{Tol}=10^{-6}$. Right: The global error for the choices $\operatorname{Tol}=10^{-3},10^{-6},10^{-9}$ (read from top to bottom). The respective overall number of steps made are 12, 77, and 856. For both pictures the initial step size was set to $h_{1,trial}=0.5$ and the parameter $\varepsilon$ was set to $1$.}
	\label{plot:Airy_WKBMM_rel_err}
\end{figure}

\begin{figure}[H]
	\centering
	\includegraphics[scale=0.95]{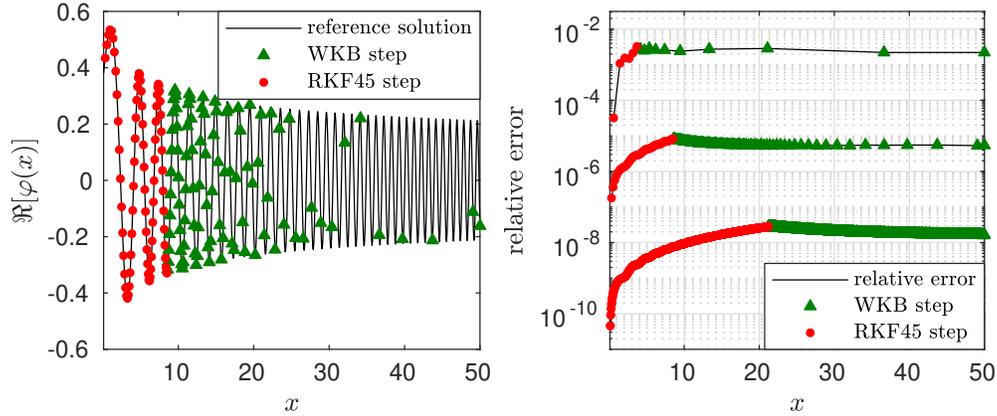}
	\caption{Left: Real part of the numerical solution obtained by using the RKWKBmod method compared to the (exact) reference solution (solid line) for $\operatorname{Tol}=10^{-6}$. Right: The global error for the choices $\operatorname{Tol}=10^{-3},10^{-6},10^{-9}$ (read from top to bottom). The respective overall number of steps made are 16, 171, and 2352. For both pictures the initial step size was set to $h_{1,trial}=0.5$, the parameter $\varepsilon$ was set to $1$ and a third order WKB-ansatz was used.}
	\label{plot:Airy_RKWKB_rel_err}
\end{figure}

\begin{figure}[H]
	\centering
	\includegraphics[scale=0.95]{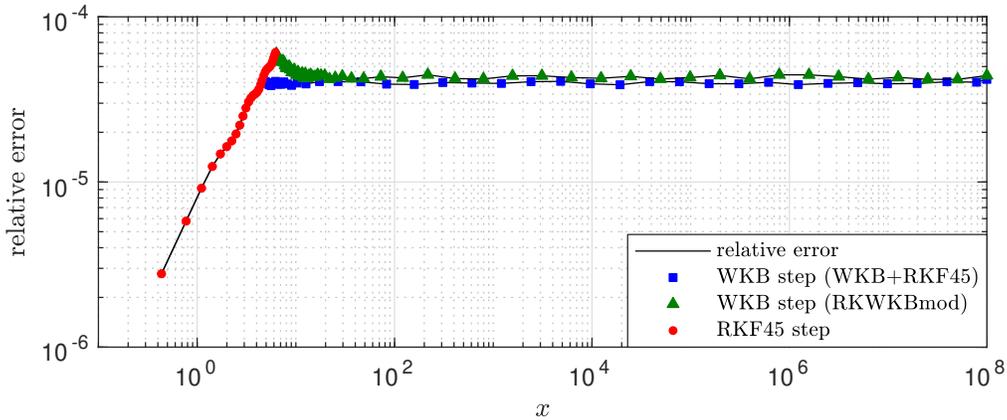}
	\caption{Global error comparison of WKB+RKF45 and RKWKBmod on $[0.1,10^{8}]$ for $\varepsilon=1$ using an initial step size $h_{1,trial}=0.5$. The error tolerance was set to $\operatorname{Tol}=10^{-5}$. A third order WKB-ansatz was used for the RKWKBmod method. Overall the algorithm using WKB+RKF45 made 58 steps, whereas RKWKBmod made 91 steps.}
	\label{plot:Airy_RKWKBMM_RKWKB_rel_err_10^8}
\end{figure}

The first example investigated in \cite{HLH16} is the Airy equation
\begin{align}\label{eqn:Airy}
	\begin{cases}
		\varepsilon^{2}\varphi^{\prime\prime  }(x) + x \varphi(x) = 0 \Comma \quad x > 0 \Comma \\
		\varphi(0) =  \frac{3^{-\frac{2}{3}} + \ii 3^{-\frac{1}{6}}}{\Gamma(\frac{2}{3})}\Comma \\
		\varepsilon\varphi^{\prime}(0) = \frac{3^{-\frac{1}{3}} - \ii 3^{\frac{1}{6}}}{\varepsilon^{-1/3}\Gamma(\frac{1}{3})} \Comma
\end{cases}
\end{align}
which results if one chooses the coefficient function $a(x)=x$. Here, $\Gamma$ denotes the gamma function. 
%\textcolor{red}{Motivation for the IC´s. Application? Beispiel eignet sich gut um die Vorteile einer WKB-Methode zu vernschaulichen, da die Lösung immer stärker oszilliert.}
The exact solution to the problem (\ref{eqn:Airy}) is given by
\begin{align}
\varphi_{exact}(x) = \Ai(-\frac{x}{\varepsilon^{2/3}}) + \ii \Bi(-\frac{x}{\varepsilon^{2/3}}) \Comma \nonumber
\end{align}
where $\Ai$ and $\Bi$ denote the Airy functions of first and second kind, respectively. This example demonstrates very well the advantages of a WKB method, since the solution becomes more and more oscillatory for large values of $x$. While standard adaptive ODE methods, e.g., Runge-Kutta methods, would need to decrease the step size more and more, a WKB-based method does not have to resolve the individual oscillations. Actually, it even allows to increase the step size for large $x$ and is therefore highly efficient. This can be seen, e.g., in Figures \ref{plot:Airy_WKBMM_rel_err}-\ref{plot:Airy_RKWKBMM_RKWKB_rel_err_10^8}.

\begin{figure}[H]
	\centering
	\includegraphics[scale=0.95]{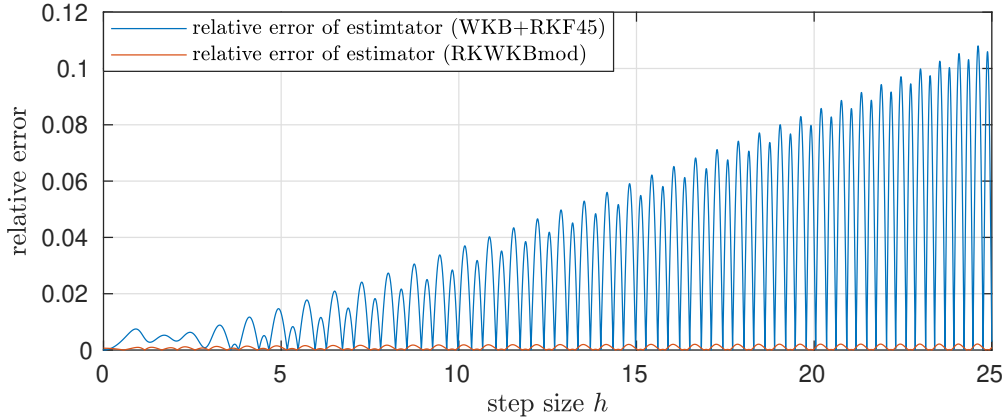}
	\caption{The relative error between the estimator (\ref{wkb_est}) and the local truncation errors for WKB+RKF45 and RKWKBmod in one single WKB step as functions of the step size $h$ for a fixed starting point $x_{0}=10$. Here, we set $\operatorname{Tol}=10^{-5}$ and $\varepsilon=1$. For the computation of the estimator with the RKWKBmod method a $2^{nd}$ and $3^{rd}$ WKB-order was used.}
	\label{plot:ErrEst_Airy_function_of_h}
\end{figure}

\begin{figure}[H]
	\centering
	\includegraphics[scale=0.95]{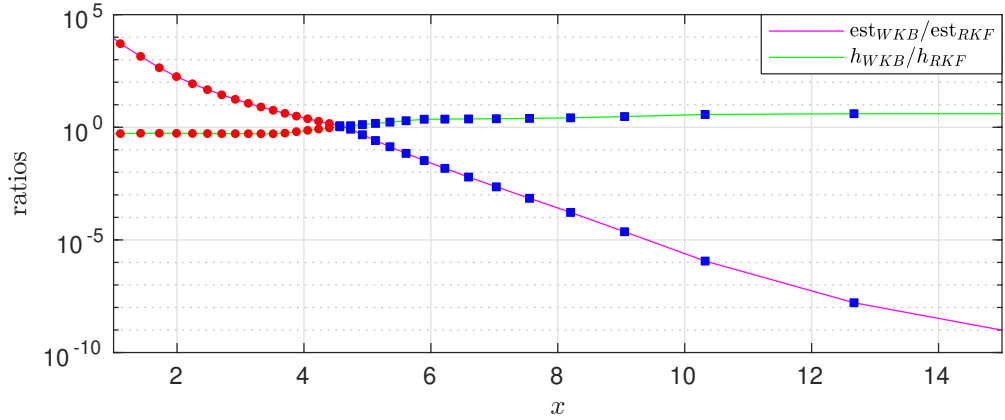}
	\caption{Numerical computation using WKB+RKF45: The magenta line indicates the ratio of the two error estimators (due to the WKB and RKF scheme). The green line gives the analogous ratio of the step sizes (locally) proposed by these two methods. For this computation we set $\operatorname{Tol}=10^{-5}$ and $\varepsilon=1$.}
	\label{plot:Airy_RKWKBMM_ratios}
\end{figure}

\begin{figure}[H]
	\centering
	\includegraphics[scale=0.95]{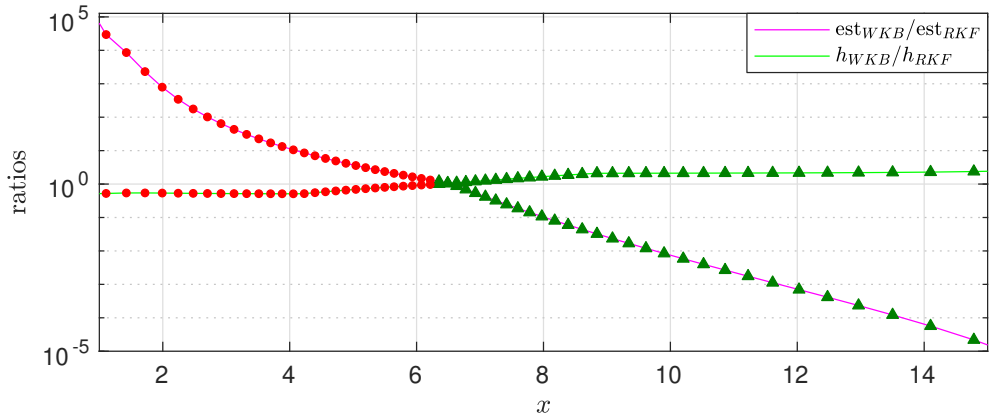}
	\caption{Numerical computation using the RKWKBmod method: The magenta line indicates the ratio of the two error estimators (due to the WKB and RKF scheme). The green line gives the analogous ratio of the step sizes (locally) proposed by these two methods. For this computation we set $\operatorname{Tol}=10^{-5}$, $\varepsilon=1$ and a third order WKB-ansatz was used.}
	\label{plot:Airy_RKWKB_ratios}
\end{figure}

\begin{figure}[H]
	\centering
	\includegraphics[scale=0.95]{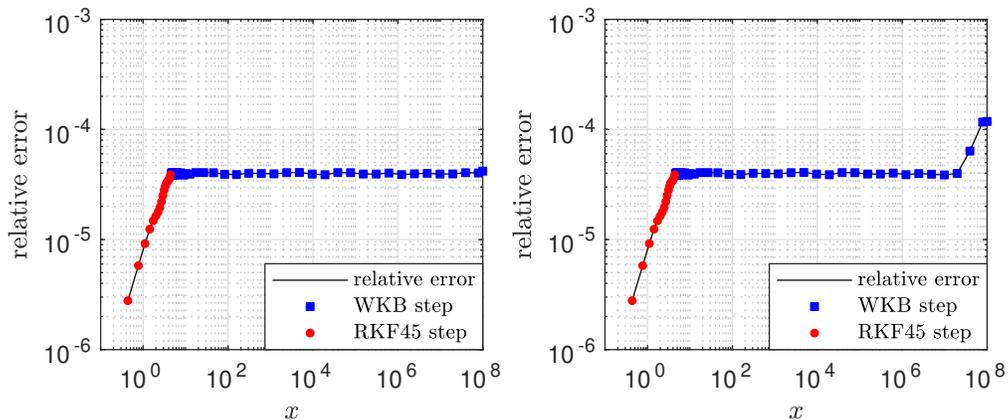}
	\caption{Left: Global (relative) error using WKB+RKF45 with exact phase. Right: Global (relative) error using WKB+RKF45 with numerically computed phase using the Clenshaw-Curtis quadrature with 15 integration nodes per step. For both computations we set $h_{1,trial}=0.5$, $\operatorname{Tol}=10^{-5}$ and $\varepsilon=1$.}
	\label{plot:Airy_WKBMM_numphase}
\end{figure}

Before starting to discuss the numerical solution of (\ref{eqn:Airy}) we first remark that the \emph{evaluation} of an oscillatory function (even of trigonometric functions such as $\sin$ or $\cos$) is numerically ill posed for very large values of $x$. This is a generic problem and not related to the numerical solution of (\ref{eqn:Airy}), or to the choice of the numerical method:
\begin{Remark}\label{remark_double_prec}
	We consider to evaluate $\varphi_{exact}(x)$ for an argument $x$ that is only known with finite accuracy, specified by machine $eps$. Then, using the lowest order expansions for $\varphi_{exact}$ and $\varphi_{exact}'$ from \ref{appendix}, we have
	$$
	\varphi_{exact}(x) \sim \frac{\varepsilon^{1/6}}{\sqrt \pi x^{1/4}}\, \exp\left(\ii\left(\frac{\pi}{4}-\frac{2}{3}\frac{x^{3/2}}{\varepsilon}\right)\right) \qquad \text{as}\quad x\to\infty \Comma
	$$
	and the relative error of $\varphi_{exact}(x(1+eps))$ is asymptotically $eps\,x^{3/2}/\varepsilon$. 
	
	Considering double precision, e.g., with \MATLAB's $eps\approx 2.2\cdot 10^{-16}$ and $x=50$, this generic evaluation error is $7.8\cdot 10^{-10}$ for $\varepsilon=10^{-4}$ and $7.8\cdot 10^{-11}$ for $\varepsilon=10^{-3}$. These unavoidable errors limit the achievable accuracy, and it will play a role in Figures \ref{plot:Airy_WKBMM_RKWKBmod_RKWKB_RKF45_tol_vs_err}-\ref{plot:Airy_WKBMM_RKWKBmod_RKWKB_RKF45_CPUtime_vs_err} below.
\end{Remark}
\medskip
As a first step we shall now test the reliability of our choice of error estimator (\ref{wkb_est}) -- but only for the WKB steps of WKB+RKF45 and RKWKBmod. Since we know the exact solution to this problem, we can compute the local truncation error in each step and are able to compare it to the error estimator. Moreover, we apply the adaptive step size control from Section \ref{chap:stepsize} with the error tolerance $\operatorname{Tol}=10^{-5}$ and set $\varepsilon=1$.
\begin{figure}[H]
	\centering
	\includegraphics[scale=0.95]{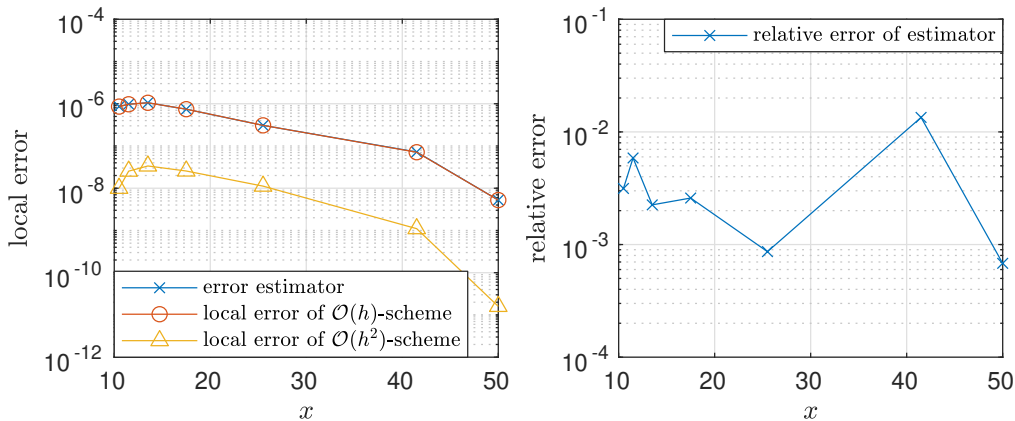}
	\caption{Left: The error estimator (\ref{wkb_est}) in comparison to the actual local truncation error computed using only WKB steps for WKB+RKF45 of first order in $h$, (\ref{scheme:first_order}), and second order in $h$, (\ref{scheme:second_order}). Right: The relative error between the estimator and the local truncation error for the $\mathcal{O}(h)$-scheme. In both pictures we set $\operatorname{Tol}=10^{-5}$ and $\varepsilon=1$.}
	\label{plot:wkb_est}
\end{figure}
\begin{figure}[H]
	\centering
	\includegraphics[scale=0.95]{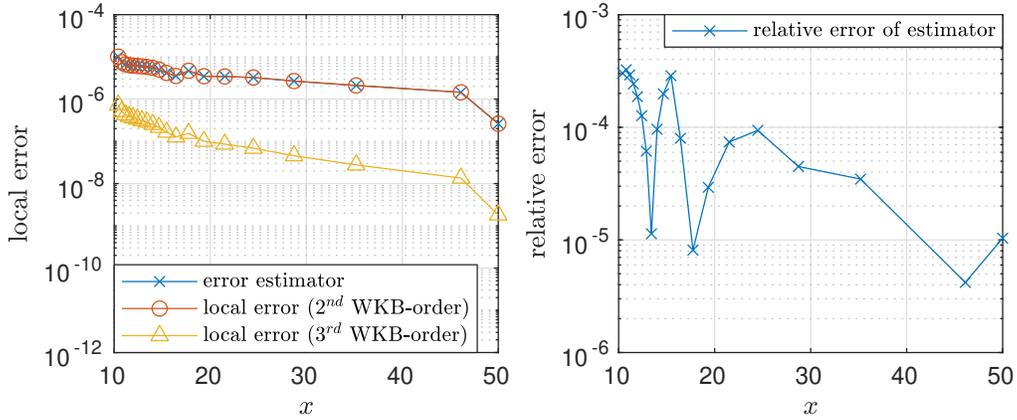}
	\caption{Left: The error estimator (\ref{wkb_est}) in comparison to the actual local truncation error computed using only WKB steps in the RKWKBmod method. Since the method is of first order (in $h$), the estimator as well as the local errors were computed by using different WKB-orders instead (here order 2 and 3). Right: The relative error between the estimator and the local truncation error ($2^{nd}$ WKB-order). In both pictures we set $\operatorname{Tol}=10^{-5}$ and $\varepsilon=1$.}
	\label{plot:rkf_est}
\end{figure}
According to the results of Figures \ref{plot:wkb_est} and \ref{plot:rkf_est}, the error estimator is in excellent agreement with the local truncation error (for this example). Hence, it seems to be an adequate choice. We also find a very good agreement by plotting the respective relative errors for one single step as functions of the step size $h$ for a fixed starting point $x_{0}$, see Figure \ref{plot:ErrEst_Airy_function_of_h}. Again, the relative error is for both methods much smaller than one. Note that in the case of WKB+RKF45 the relative error goes to zero, if $h$ tends to zero. Therefore, the estimator seems to be asymptotically correct.
\begin{Remark}
	The computation of the local truncation error involves the evaluation of the Airy function, which seems to have an accuracy problem for large values of $x$ when using the standard routine \verb+airy()+ in \MATLAB. Hence we used a modified implementation of the Airy functions, which consists of the function \verb+airy()+ from \MATLAB for small values of $x$ and asymptotic expressions for the Airy function for large values of $x$. More precisely, the evaluations were performed as given in Table \ref{table:xsp}. The asymptotic expressions are based on well-known expansions, which can be found in \ref{appendix}. The order for truncating the expansions was set to $K=3$.
	\begin{table}[H] 
		\centering
		\begin{tabular}{|l|l|l|}
			\hline
			&     \verb+airy()+  & asymptotics   \\ \hline
			$\Ai$       &   $x\in[0,500]$    & $x\in(500,\infty)$\\ \hline
			$\Ai'$      & 	$x\in[0,400]$    & $x\in(400,\infty)$\\ \hline
			$\Bi$       &   $x\in[0,500]$    & $x\in(500,\infty)$\\ \hline
			$\Bi'$      &   $x\in[0,400]$    & $x\in(400,\infty)$\\ \hline
		\end{tabular}
		
		\caption{Intervals for evaluating the Airy function of first and second kind as well as their derivatives. For small values of $x$ the original function from \MATLAB was used, and for large $x$ the evaluation was performed using the asymptotic expressions (\ref{Asymptotic_Ai})-(\ref{Asymptotic_Bi'}) with truncation after $K=3$.}
		
		\label{table:xsp}
	\end{table}
\end{Remark}
Now, we will give numerical results for solving the initial value problem (\ref{eqn:Airy}) for the Airy equation with WKB+RKF45 and RKWKBmod. We recall that the algorithm automatically chooses between RKF45 and the respective WKB steps. To illustrate this difference in the Figures \ref{plot:Airy_WKBMM_rel_err}-\ref{plot:Airy_WKBMM_numphase} as well as Figures \ref{plot:PCF_WKBMM_eps=2^-6}-\ref{plot:PCF_RKWKBmod_eps=2^-6}, we will mark RKF45 steps with red dots, WKB steps using WKB+RKF45 with blue squares, and WKB steps using the RKWKBmod method with green triangles. To exclude the turning point $x=0$, we solve the Airy equation (\ref{eqn:Airy}) in Figures \ref{plot:Airy_WKBMM_rel_err} and \ref{plot:Airy_RKWKB_rel_err} on $[0.1, 50]$.
%\part{\begin{figure}[H]
%	\centering
%	\subfloat[RKWKB: 42 steps]{{\includegraphics[scale=.7]{x_start=1_to_10^8_RKWKB_alpha=10^-4.png} }}%
%	%\includegraphics[scale=0.7]{x_start=1_to_10^8_RKWKB_alpha=10^-4.png}
%	\qquad
%	\subfloat[RKWKBMM: 25 steps]{{\includegraphics[scale=.7]{x_start=1_to_10^8_RKWKBMM_alpha=10^-4.png} }}%
%	%\includegraphics[scale=0.7]{x_start=1_to_10^8_RKWKBMM_alpha=10^-4.png}
%	
%	\caption{(a): The algorithm switch at $x\approx 3.2$ from the RKF4(5) scheme to WKB. The relative error starts at around $\sim 10^{-6}$ on the RK side and goes up to $\sim 2\cdot 10^{-4}$ on the WKB side. (b) The algorithm switches at $x\approx 2.5$ from RK to WKB. The error only goes up to $\sim 2\cdot 10^{-5}$ on the WKB side, so it is better by a factor of at least $10$.}
%	\label{img:grafik-dummy}
%\end{figure}}
According to Figures \ref{plot:Airy_WKBMM_rel_err} and \ref{plot:Airy_RKWKB_rel_err}, WKB+RKF45 seems to perform slightly better than the RKWKBmod method in this example, in matters of global error. But at the same time WKB+RKF45 needs significantly fewer steps than RKWKBmod. Also, within the algorithm using WKB+RKF45, the switch from RKF45 steps to WKB steps happens earlier as can be seen in Figure \ref{plot:Airy_WKBMM_rel_err}. In Figure \ref{plot:Airy_RKWKBMM_RKWKB_rel_err_10^8}, the relative global errors of both algorithms are compared on $[0.1,10^{8}]$. We find that the algorithm using WKB+RKF45 made fewer steps (58 vs.\ 91) while producing a slightly lower global error at the same time. The almost identical grid spacing of both methods from $x\approx 300$ onwards is due to limiting the quotient of two consecutive step sizes, imposed in both methods.
Furthermore, Figures \ref{plot:Airy_RKWKBMM_ratios} and \ref{plot:Airy_RKWKB_ratios} show the ratio between the WKB- and RK- error estimators as well as the ratio of the proposed step sizes (by the WKB and RKF schemes) for each step of the above simulations.
For both WKB+RKF45 and RKWKBmod, these plots demonstrate well the superiority of the RK scheme in the less oscillatory regime, i.e.\ for $x$ small, as the ratio of the error estimators is very large there. The switching points to the WKB schemes are well defined (again for both WKB methods) -- due to the monotonous behaviour of both the ratio of error estimators and the ratio of proposed step sizes. Note also that the step size ratios are bounded from above and below because of the limitations in (\ref{theta}).

We also want to give a perception of how WKB+RKF45 can perform, if the phase integral (\ref{phase}) can not be evaluated exactly. For this purpose we use the \textit{Clenshaw-Curtis} quadrature (cf.\ \cite{CC60}) to approximate the phase (\ref{phase}). That spectral method was already used in combination with the WKB-marching method, see \cite[§5]{AKU19}. In \cite{AKU19} they approximate the phase at first on Chebyshev collocation points throughout the \textit{whole interval} and use barycentric interpolation for the ODE-grid points $x_{n}$. This was possible since they worked only on the ``small'' interval $[0,1]$. In contrast, we will instead approximate the phase (\ref{phase}) in \textit{each interval} $[x_{n},x_{n+1}]$ individually, since we are dealing here with the ``long'' interval $[0.1,10^8]$. Figure \ref{plot:Airy_WKBMM_numphase} gives a comparison of the global error for the Airy equation (\ref{eqn:Airy}), when using the exact phase vs.\ a numerically computed phase. According to these results the relative $\varphi$-errors from computing the phase (\ref{phase}) numerically become visible only from $x\approx 10^{7}$ onwards.

\medskip
We will next compare the numerical results of four methods, namely, WKB+RKF45, RKWKBmod, RKWKB, and RKF45 on the Airy equation for several values of the parameter $\varepsilon$. For RKF45, we do not use any smaller value of $\varepsilon$ than $10^{-1}$, as the CPU time gets exorbitantly large. In Figure \ref{plot:Airy_WKBMM_RKWKBmod_RKWKB_RKF45_tol_vs_err} we compare the accuracy of each of the methods depending on $\operatorname{Tol}$. We find that WKB+RKF45 outperforms the other methods, in matters of global errors, particularly for small values of $\varepsilon$ and $\operatorname{Tol}$. By Remark \ref{remark_double_prec} the accuracy limit at $x=50$ is $7.8\cdot 10^{-10}$ for $\varepsilon=10^{-4}$, which seems to explain the lower bound of the errors in Figure \ref{plot:Airy_WKBMM_RKWKBmod_RKWKB_RKF45_tol_vs_err}. Note also that for RKF45 the error increases like $\mathcal{O}(\varepsilon^{-1})$ for smaller $\varepsilon$-values, while it decreases for the WKB methods.

In Figure \ref{plot:Airy_WKBMM_RKWKBmod_RKWKB_RKF45_CPUtime_vs_err} we give a work-precision diagram; for a fair comparison between the methods, points showing the same error should be compared. There is a big difference between WKB+RKF45 and RKWKB(mod), regarding the CPU time: For $\varepsilon=10^{0},10^{-1},10^{-2}$ this difference is particularly significant for small errors (stemming from small prescribed tolerances). For $\varepsilon=10^{-3}$ WKB+RKF45 already beats RKWKB(mod) for all data points. For $\varepsilon=10^{-4}$ the error intervals of WKB+RKF45 and RKWKB(mod) do not overlap. But for the same CPU time WKB+RKF45 yields much more accurate results than RKWKB or RKWKBmod. Overall we conclude from Figure \ref{plot:Airy_WKBMM_RKWKBmod_RKWKB_RKF45_CPUtime_vs_err} that WKB+RKF45 outperforms RKWKBmod and RKWKB significantly. Using RKF45, the CPU time increases drastically for smaller $\varepsilon$-values.

\begin{figure}[H]
	\centering
	\includegraphics[scale=0.95]{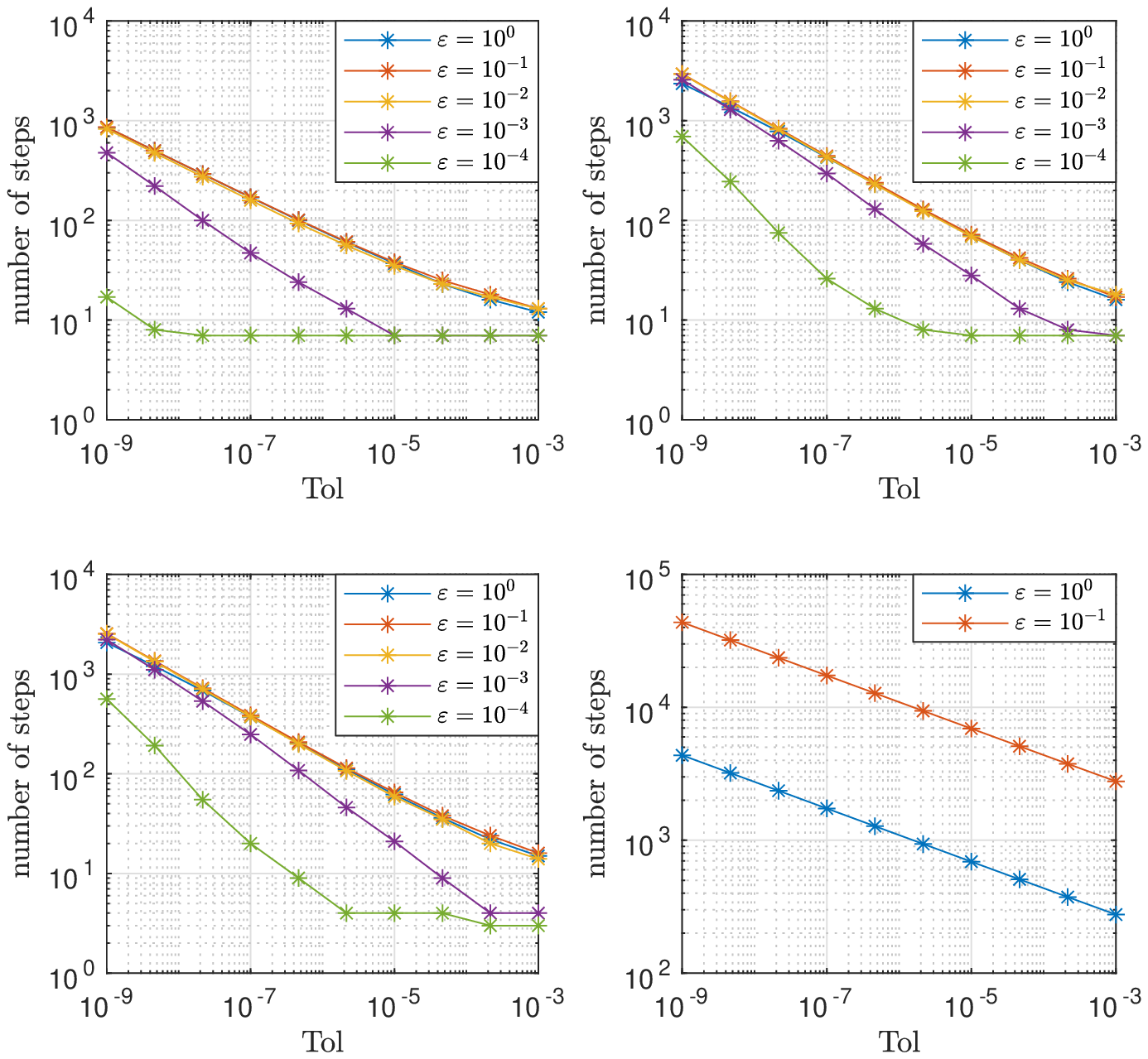}
	\caption{Number of steps needed for the Airy equation (\ref{eqn:Airy}) on $[0.1, 50]$ as a function of $\operatorname{Tol}$, for several $\varepsilon$-values. Top-left: WKB+RKF45. Top-right: RKWKBmod. Bottom-left: RKWKB. Bottom-right: RKF45. For all methods we set $h_{1,trial}=0.5$. A third order WKB-ansatz for RKWKBmod and RKWKB was used.}
	\label{plot:Airy_WKBMM_RKWKBmod_RKWKB_RKF45_tol_vs_steps}
\end{figure}

Figure \ref{plot:Airy_WKBMM_RKWKBmod_RKWKB_RKF45_tol_vs_steps} displays the respective number of steps needed in each method, as a function of the prescribed tolerance. Note that for WKB+RKF45 and RKWKBmod the number of steps is bounded from below. This is because of limiting the quotient of two consecutive step sizes and it is clearly visible for $\varepsilon=10^{-3},10^{-4}$.
\begin{figure}[H]
	\centering
	\includegraphics[scale=0.95]{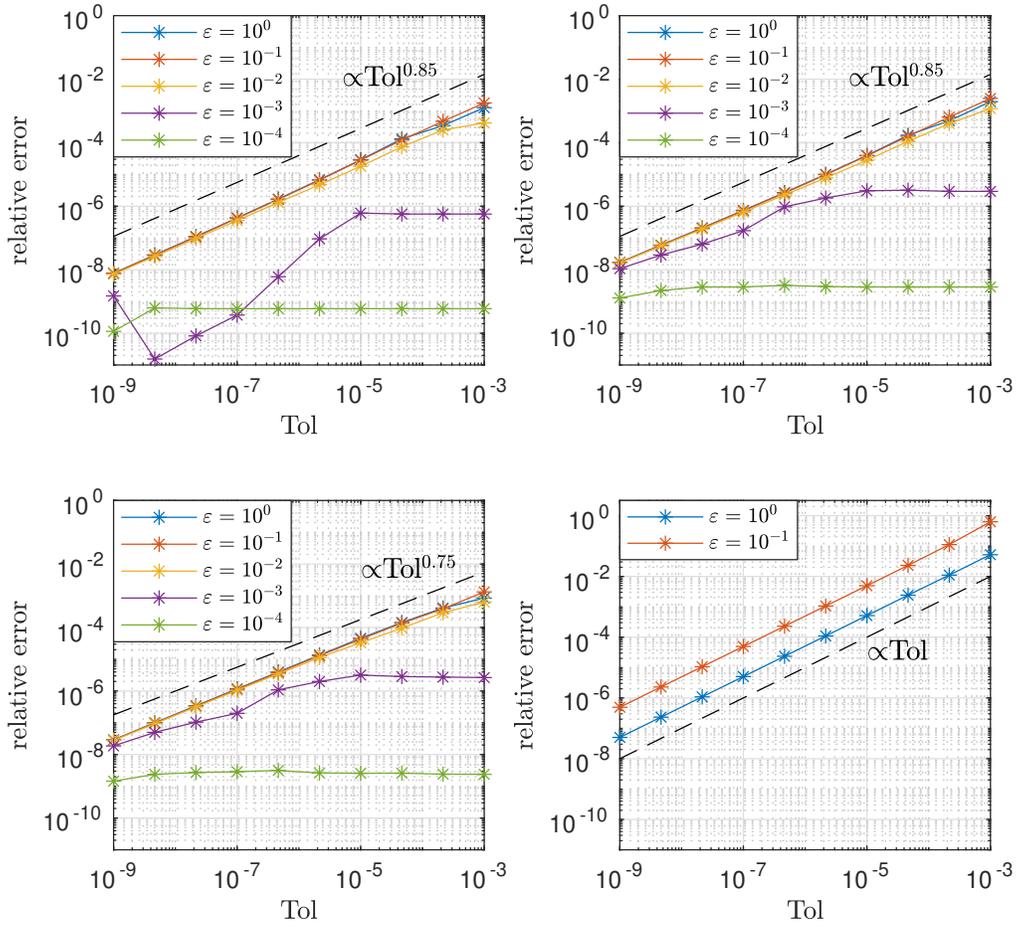}
	\caption{Global (relative) errors (in the $l^{2}$-norm) for the Airy equation (\ref{eqn:Airy}) on $[0.1, 50]$ as a function of $\operatorname{Tol}$, for several $\varepsilon$-values. Top-left: WKB+RKF45. Top-right: RKWKBmod. Bottom-left: RKWKB. Bottom-right: RKF45. For all methods we set $h_{1,trial}=0.5$. A third order WKB-ansatz for RKWKBmod and RKWKB was used.}
	\label{plot:Airy_WKBMM_RKWKBmod_RKWKB_RKF45_tol_vs_err}
\end{figure}
\begin{figure}[H]
	\centering
	\includegraphics[scale=0.95]{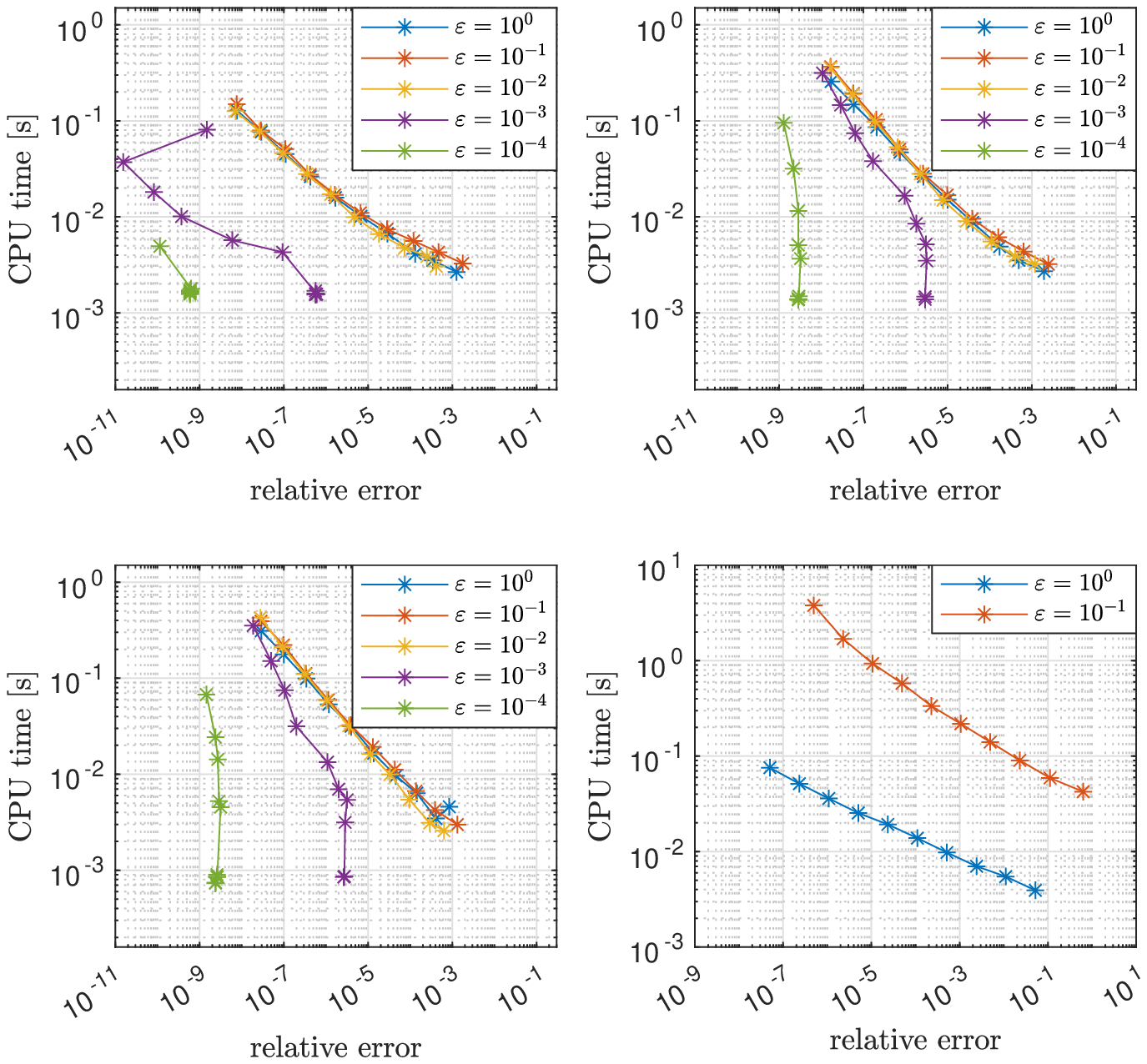}
	\caption{CPU times vs.\ global (relative) errors (in the $l^{2}$-norm) for the Airy equation (\ref{eqn:Airy}) on $[0.1, 50]$, computed for 10 logarithmically evenly spaced values of $\operatorname{Tol}$ in the range $[10^{-9},10^{-3}]$, for several $\varepsilon$-values. Top-left: WKB+RKF45. Top-right: RKWKBmod. Bottom-left: RKWKB. Bottom-right: RKF45 (note the different scales). For all methods we set $h_{1,trial}=0.5$. A third order WKB-ansatz for RKWKBmod and RKWKB was used.}
	\label{plot:Airy_WKBMM_RKWKBmod_RKWKB_RKF45_CPUtime_vs_err}
\end{figure}
\subsection{Second Example: Parabolic cylinder function}
\begin{figure}[H]
	\centering
	\includegraphics[scale=0.95]{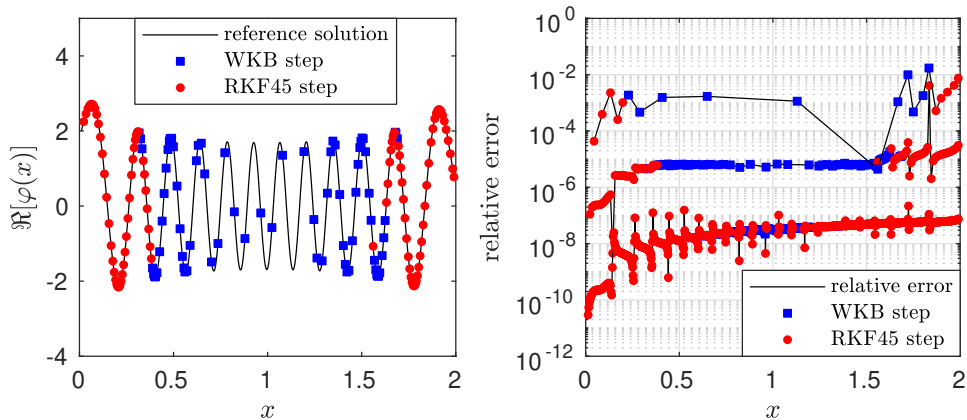}
	\caption{Left: Real part of the numerical solution obtained by using WKB+RKF45 compared to the (exact) reference solution (solid line) for $\operatorname{Tol}=10^{-6}$. Right: The global error for the choices $\operatorname{Tol}=10^{-3},10^{-6},10^{-9}$ (read from top to bottom). The respective overall number of steps made are 21, 166, and 1287. For both pictures the initial step size was set to $h_{1,trial}=0.05$ and the parameter $\varepsilon$ was set to $2^{-6}$.}
	\label{plot:PCF_WKBMM_eps=2^-6}
\end{figure}

\begin{figure}[H]
	\centering
	\includegraphics[scale=0.95]{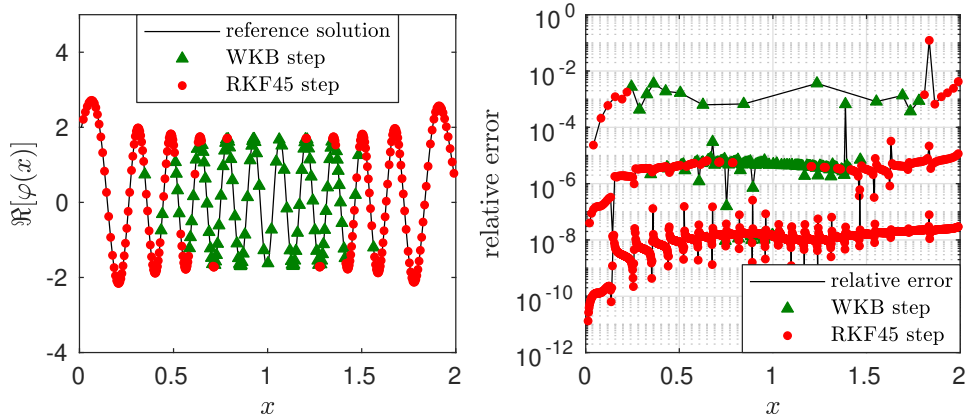}
	\caption{Left: Real part of the numerical solution obtained by using RKWKBmod compared to the (exact) reference solution (solid line) for $\operatorname{Tol}=10^{-6}$. Right: The global error for the choices $\operatorname{Tol}=10^{-3},10^{-6},10^{-9}$ (read from top to bottom). The respective overall number of steps made are 26, 326, and 1543. For both pictures the initial step size was set to $h_{1,trial}=0.05$, the parameter $\varepsilon$ was set to $2^{-6}$ and a third order WKB-ansatz was used.}
	\label{plot:PCF_RKWKBmod_eps=2^-6}
\end{figure}

The second example is taken from \cite{AD20} and includes a quadratic coefficient function $a(x)$:
%\textcolor{blue}{($\kappa$ ist frei wählbar und dient hier nur eine Runterskalierung, hier ist $a(x)=k_{1}x^{2}+k_{2}x$ mit $k_{1}=-\frac{1}{2}$ und $k_{2}=1$.)}
%
\begin{align}\label{eqn:pcf}
	\begin{cases}
		\varepsilon^{2}\varphi^{\prime\prime}(x) + \left(-\frac{1}{2}x^{2}+x\right) \varphi(x) = 0 \Comma \quad x\in (0,2) \Comma\\
		\varphi(0) = \kappa U(\nu,z(0)) \Comma \\
		\varphi^{\prime}(0) = -\kappa 2^{-\frac{1}{4}}\varepsilon^{-\frac{1}{2}}U'(\nu,z(0)) \Comma
	\end{cases}
\end{align}
with
\begin{align}
	\nu:=-\frac{1}{\sqrt{8}\varepsilon} \Comma \quad z(x):=\frac{2^{\frac{1}{4}}}{\sqrt{\varepsilon}}(1-x) \Comma \nonumber
\end{align}
and
\begin{align}
	\kappa :=\frac{2}{U(\nu,0)-\ii\sqrt{\varepsilon}2^{\frac{3}{4}}U^{\prime}(\nu,0)} \period \nonumber
\end{align}
The exact solution reads
\begin{align}
	\varphi_{exact}(x)=\kappa U(\nu,z(x)) \Comma \nonumber
\end{align}
where $U(\nu,z)$ denotes the parabolic cylinder function (PCF) (cf.\ \cite[§12]{NHM10}). As before, let us compare numerical results only for WKB+RKF45 and RKWKBmod at first. There are two turning points, namely at $x=0$ and $x=2$. Therefore, we expect the two methods to make RKF45 steps near the turning points and WKB steps between them.
%The Figures \ref{plot:PCF_RKWKBMM_eps=2^-6}-\ref{plot:PCF_RKWKB_eps=2^-6} demonstrate how the parameter $\varepsilon$ affects the solution. As $\varepsilon$ gets smaller the solution becomes more and more oscillating.
Numerical results for the specific choice $\varepsilon=2^{-6}$ are presented in Figures \ref{plot:PCF_WKBMM_eps=2^-6}-\ref{plot:PCF_RKWKBmod_eps=2^-6}.
According to Figures \ref{plot:PCF_WKBMM_eps=2^-6} and \ref{plot:PCF_RKWKBmod_eps=2^-6} no significant difference between WKB+RKF45 and RKWKBmod can be observed for $\varepsilon=2^{-6}$, in matters of global error. But again, WKB+RKF45 needs significantly fewer steps than RKWKBmod. Within the algorithm using WKB+RKF45, the switch-over between RKF45 steps and WKB steps happens closer to the turning points.

\begin{figure}[H]
	\centering
	\includegraphics[scale=0.95]{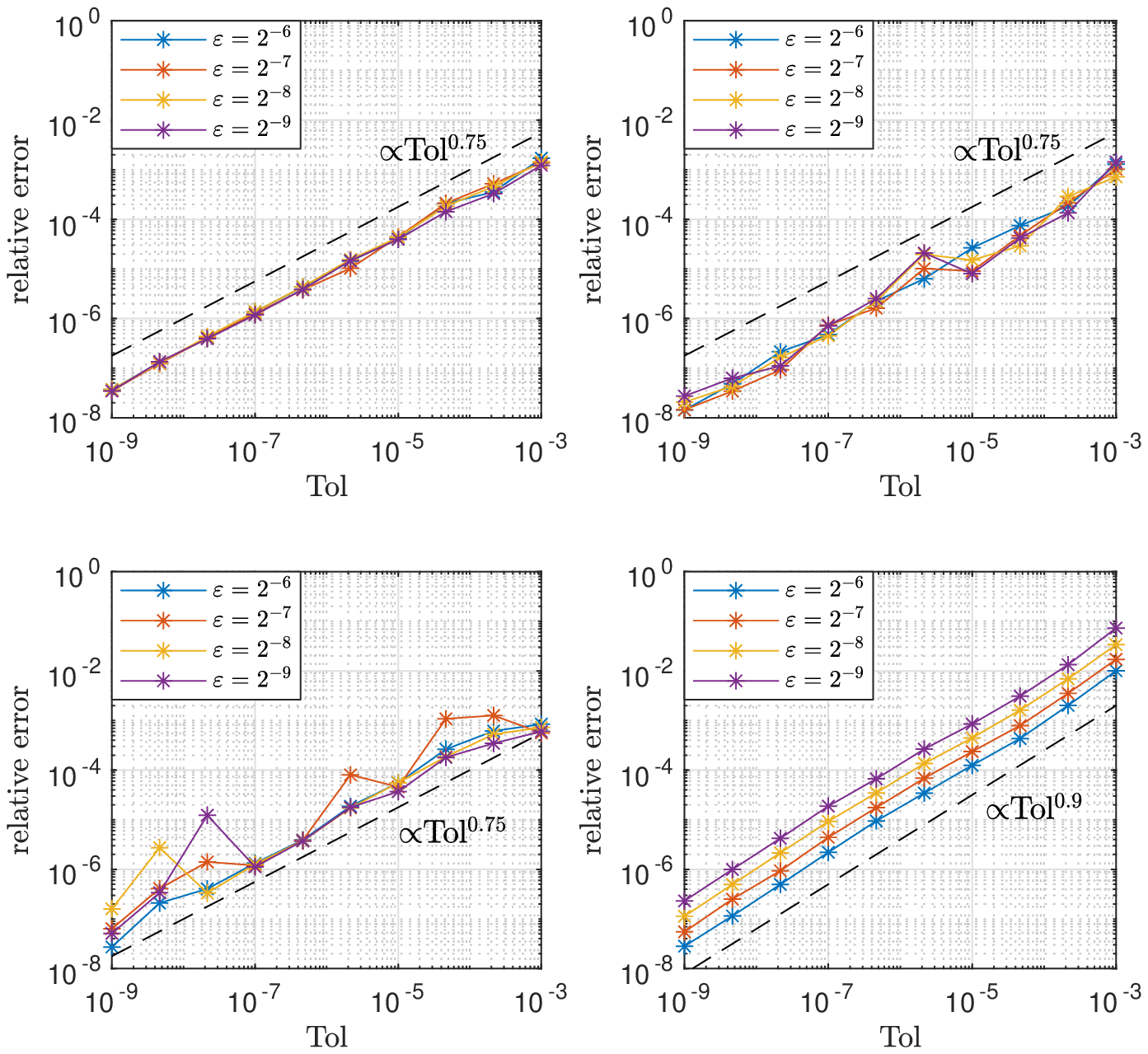}
	\caption{Global (relative) errors (in the $l^{2}$-norm) for equation (\ref{eqn:pcf}) on $[0.01, 1.99]$ as a function of $\operatorname{Tol}$, for several $\varepsilon$-values. Top-left: WKB+RKF45. Top-right: RKWKBmod. Bottom-left: RKWKB. Bottom-right: RKF45. For all methods we set $h_{1,trial}=0.05$. A third order WKB-ansatz for RKWKBmod and RKWKB was used.}
	\label{plot:PCF_WKBMM_RKWKBmod_RKWKB_RKF45_tol_vs_err}
\end{figure}

We shall now present numerical results for the four methods WKB+RKF45, RKWKBmod, RKWKB, and RKF45 for several values of $\varepsilon$. In Figure \ref{plot:PCF_WKBMM_RKWKBmod_RKWKB_RKF45_tol_vs_err}, we compare the global errors of each method depending on $\operatorname{Tol}$. Here, we observe only a small difference between WKB+RKF45 and RKWKBmod. Using RKWKB, less smooth error curves can bee seen, with large peaks for lower errors (i.e.\ lower tolerances). In contrast to every other method, WKB+RKF45 seems to produce quite $\varepsilon$-independent global error curves. For RKF45 one sees that the error again increases like $\mathcal{O}(\varepsilon^{-1})$ for smaller values of $\varepsilon$.
	
\begin{figure}[H]
	\centering
	\includegraphics[scale=0.95]{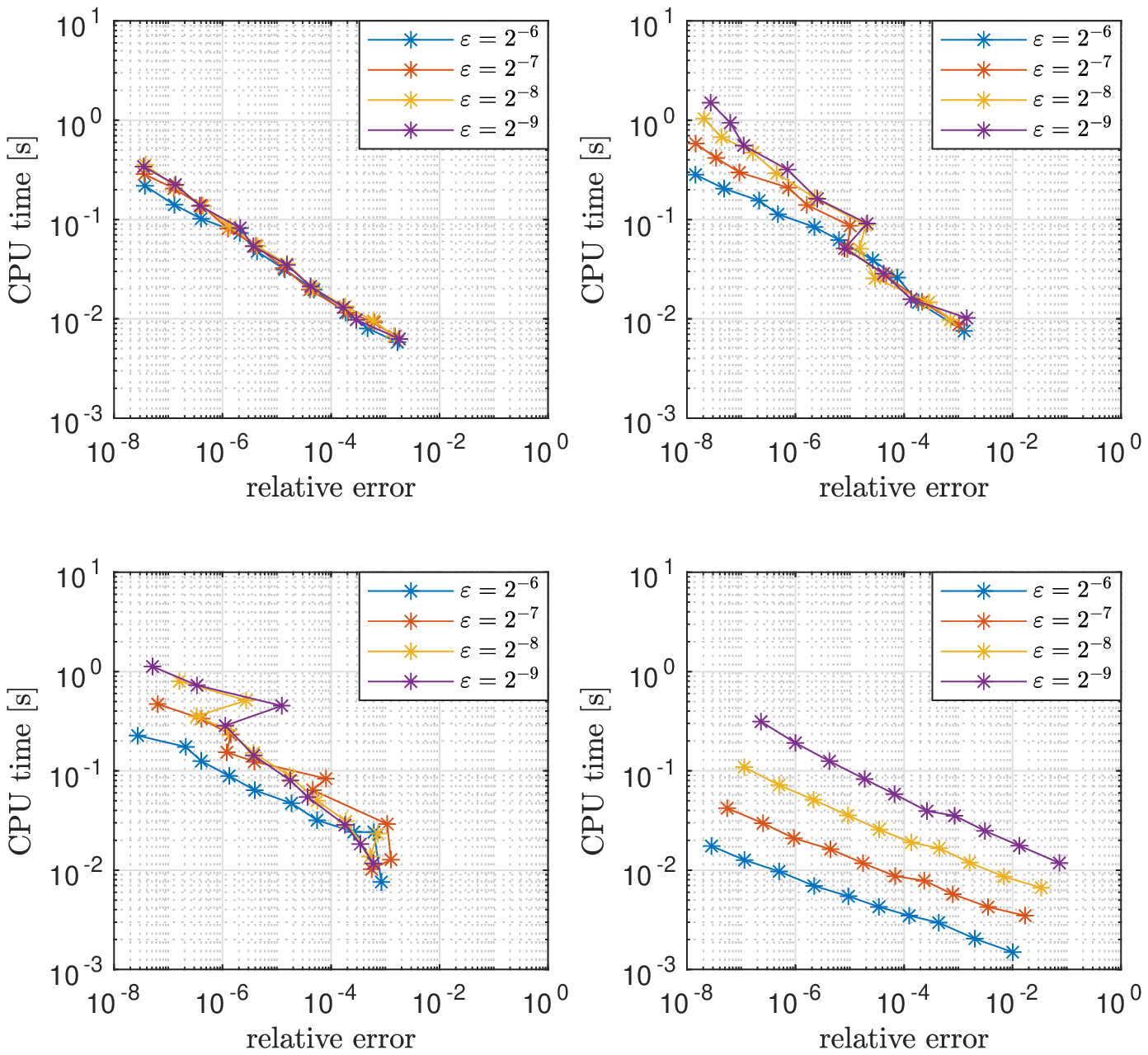}
	\caption{CPU times vs.\ global (relative) errors (in the $l^{2}$-norm) for equation (\ref{eqn:pcf}) on $[0.01, 1.99]$, computed for 10 logarithmically evenly spaced values of $\operatorname{Tol}$ in the range $[10^{-9},10^{-3}]$, for several $\varepsilon$-values. Top-left: WKB+RKF45. Top-right: RKWKBmod. Bottom-left: RKWKB. Bottom-right: RKF45. For all methods we set $h_{1,trial}=0.05$. A third order WKB-ansatz for RKWKBmod and RKWKB was used.}
	\label{plot:PCF_WKBMM_RKWKBmod_RKWKB_RKF45_CPUtime_vs_err}
\end{figure}

In the work-precision diagrams in Figure \ref{plot:PCF_WKBMM_RKWKBmod_RKWKB_RKF45_CPUtime_vs_err}, we observe for WKB+RKF45 that the CPU times are quite independent of $\varepsilon$, whereas for RKWKB(mod) they grow with decreasing $\varepsilon$, particularly for small errors (stemming from small prescribed tolerances). Overall we conclude from Figure \ref{plot:PCF_WKBMM_RKWKBmod_RKWKB_RKF45_CPUtime_vs_err} that WKB+RKF45 outperforms RKWKBmod and RKWKB (particularly for small $\varepsilon$ and small tolerances) while showing an $\varepsilon$-uniform behaviour at the same time. Note also that, for RKF45, the CPU time increases drastically for smaller values of $\varepsilon$.

\begin{figure}[H]
	\centering
	\includegraphics[scale=0.95]{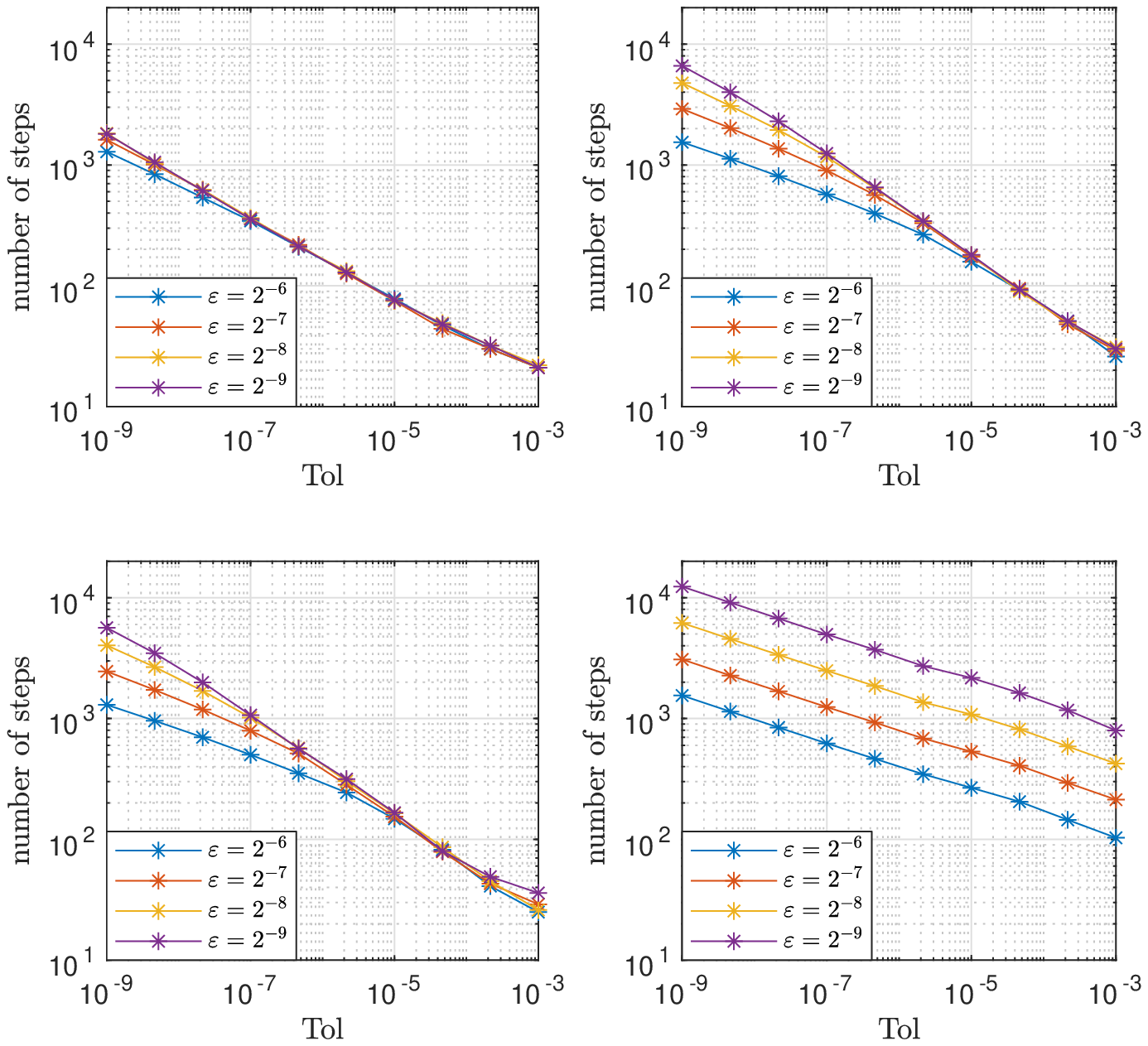}
	\caption{Number of steps needed for equation (\ref{eqn:pcf}) on $[0.01, 1.99]$ as a function of $\operatorname{Tol}$, for several $\varepsilon$-values. Top-left: WKB+RKF45. Top-right: RKWKBmod. Bottom-left: RKWKB. Bottom-right: RKF45. For all methods we set $h_{1,trial}=0.05$. A third order WKB-ansatz for RKWKBmod and RKWKB was used.}
	\label{plot:PCF_WKBMM_RKWKBmod_RKWKB_RKF45_tol_vs_steps}
\end{figure}

Figure \ref{plot:PCF_WKBMM_RKWKBmod_RKWKB_RKF45_tol_vs_steps} shows the respective number of steps needed in each method, as a function of the prescribed tolerance.
\begin{Remark}
	The computation of the reference solution involves the evaluation of the PCF, which is not readily available in \MATLAB. But the PCF can be related to the Kummer confluent hypergeometric function ${}_{1}F_{1}$ (see \cite[§13]{NHM10}), which is available in \MATLAB as \verb+hypergeom()+. However, for small parameters $\varepsilon$ this evaluation is very time consuming. For Figures \ref{plot:PCF_WKBMM_RKWKBmod_RKWKB_RKF45_tol_vs_err}-\ref{plot:PCF_WKBMM_RKWKBmod_RKWKB_RKF45_tol_vs_steps} we thus computed the reference solutions by solving the initial value problem (\ref{eqn:pcf}) with \MATLAB's routine \verb+ode45()+ (for all $\varepsilon$) using a very small error tolerance.
\end{Remark}
\section{Conclusion}\label{chap:conclusion}
We have introduced in this paper an extension to the WKB-marching method from \cite{ABN11} by including into the algorithm an adaptive step size controller as well as a switching mechanism. In numerical simulations based on two examples this method yielded smaller global errors (particularly for small tolerances and small $\varepsilon$-values) in comparison to the Runge-Kutta-WKB method from \cite{HLH16,AHLH20}, an alternative WKB-based scheme. Our tests revealed that the efficiency gain is mostly due to the different WKB method used here, while the different step size controls (here vs.\ \cite{HLH16,AHLH20}) do not play a big role. Our switching mechanism ensures well defined switching points between WKB steps and Runge-Kutta-Fehlberg 4(5) steps for oscillatory and, respectively, smoother regions of the ODE-solution. Especially for the Airy equation on the large spatial interval $[0.1,10^{8}]$ the efficiency of the method is demonstrated very well, as the scheme skips millions of oscillations within one step, while staying accurate at the same time. There is also a \MATLAB program available in a GitHub repository\footnote[1]{\href{https://github.com/JannisKoerner/adaptive-WKB-marching-method}{https://github.com/JannisKoerner/adaptive-WKB-marching-method}}, which also offers the possibility to compute the phase (\ref{phase}) numerically, as done for Figure \ref{plot:Airy_WKBMM_numphase}.
%
%\section*{CRediT authorship contribution statement \textcolor{blue}{(optional)}}
%\textbf{Jannis Körner:} blabla. \textbf{Anton Arnold:} blabla. \textbf{Kirian Döpfner:} blabla. \textcolor{blue}{(Möglich sind hier die Ausdrücke: \textbf{Conceptualization, Methodology, Software, Validation, Formal analysis, Investigation, Resources, Data Curation, Writing - Original Draft, Writing - Review \& Editing, Visualization, Supervision, Project administration, Funding acquisition}, Definitionen der Ausdrücke unter \\ https://www.elsevier.com/authors/policies-and-guidelines/credit-author-statement)}
%
%\section*{Declaration of competing interest}
%
%The authors declare that they have no known competing financial interests or personal relationships that could have appeared to influence the work reported in this paper.
%
\section*{Acknowledgements}
The authors A. Arnold and K. Döpfner were partially supported by the binational FWF-project I3538-N32. Moreover, the authors J. Körner and K. Döpfner have been (partially) supported from the Austrian Science Fund (FWF) through grant number W1245.
\appendix
\section{Asymptotic formulas for Airy functions}
\label{appendix}
For real-valued $x$ and $x \to \infty$, asymptotic expansions for the Airy functions and their first derivatives are given in \cite[§9.7 (ii)]{NHM10}:
\begin{align}
	\Ai(-x) &\sim \frac{1}{\sqrt{\pi}x^{\frac{1}{4}}}\left( \cos\left( \zeta - \frac{\pi}{4} \right) \sum_{k = 0}^{\infty} (-1)^{k} \frac{u_{2k}}{\zeta^{2k}} + \sin\left( \zeta - \frac{\pi}{4} \right) \sum_{k = 0}^{\infty} (-1)^{k} \frac{u_{2k+1}}{\zeta^{2k+1}}\right) \Comma  \label{Asymptotic_Ai}\\
	\Ai^{\prime}(-x) &\sim \frac{x^{\frac{1}{4}}}{\sqrt{\pi}}\left( \sin\left( \zeta - \frac{\pi}{4} \right) \sum_{k = 0}^{\infty} (-1)^{k} \frac{v_{2k}}{\zeta^{2k}} - \cos\left( \zeta - \frac{\pi}{4} \right) \sum_{k = 0}^{\infty} (-1)^{k} \frac{v_{2k+1}}{\zeta^{2k+1}}\right) \Comma \label{Asymptotic_Ai'} \\
	\Bi(-x) &\sim \frac{1}{\sqrt{\pi}x^{\frac{1}{4}}}\left( -\sin\left( \zeta - \frac{\pi}{4} \right) \sum_{k = 0}^{\infty} (-1)^{k} \frac{u_{2k}}{\zeta^{2k}} + \cos\left( \zeta - \frac{\pi}{4} \right) \sum_{k = 0}^{\infty} (-1)^{k} \frac{u_{2k+1}}{\zeta^{2k+1}}\right) \Comma \label{Asymptotic_Bi} \\
	\Bi^{\prime}(-x) &\sim \frac{x^{\frac{1}{4}}}{\sqrt{\pi}}\left( \cos\left( \zeta - \frac{\pi}{4} \right) \sum_{k = 0}^{\infty} (-1)^{k} \frac{v_{2k}}{\zeta^{2k}} + \sin\left( \zeta - \frac{\pi}{4} \right) \sum_{k = 0}^{\infty} (-1)^{k} \frac{v_{2k+1}}{\zeta^{2k+1}}\right) \period \label{Asymptotic_Bi'}
\end{align}
Here, $\zeta := \frac{2}{3} x^{\frac{3}{2}}$ and the coefficients $u_{k}$ and $v_{k}$ are given by (see \cite[§9.7 (i)]{NHM10}):
\begin{align}
	u_{0}&=v_{0}=1 \Comma \nonumber \\
	u_{k} &= \frac{(2 k + 1) \cdot (2 k + 3) \cdot (2 k + 5) \cdot ... \cdot (6 k - 1)}{216^{k} k!} \nonumber \\ &= \frac{(6 k - 5) \cdot (6 k - 3) \cdot (6 k - 1)}{(2 k -1) 216 k} u_{k - 1} \Comma \quad k=1,2,...\Comma\nonumber \\
	v_{k} &= \frac{6 k + 1}{1 - 6 k}u_{k}\Comma \quad k=1,2,...\period \nonumber
\end{align}
%
%\subsection{Relations of the Parabolic Cylinder Function}
%\textcolor{blue}{Ich bin mir unsicher bezüglich der Notwendigkeit dieses Unterabschnitts...}
%The PCF can be expressed in terms of the so-called Kummer confluent hypergeometric function ${}_{1}F_{1}$, see \cite[§13]{NHM10}. For real or complex valued $\nu$ and complex valued $z$ it holds
%\begin{align}
%	U(\nu,z)=U(\nu,0)u_{1}(\nu,z)+U'(\nu,0)u_{2}(\nu,z)\comma
%\end{align}
%%
%where
%%
%\begin{align}
%	u_{1}(\nu,z)&=\operatorname{e}^{-\frac{1}{4}z^{2}}{}_{1}F_{1}(\frac{1}{2}\nu+\frac{1}{4},\frac{1}{2},\frac{1}{2}z^{2})\comma\nonumber\\
%	u_{2}(\nu,z)&=z\operatorname{e}^{-\frac{1}{4}z^{2}}{}_{1}F_{1}(\frac{1}{2}\nu+\frac{3}{4},\frac{3}{2},\frac{1}{2}z^{2})\comma\nonumber\\
%	U(\nu,0)&=\frac{\sqrt{\pi}}{2^{\frac{1}{2}\nu+\frac{1}{4}}\Gamma\left(\frac{3}{4}+\frac{1}{2}\nu\right)}\comma\nonumber\\
%	U'(\nu,0)&=\frac{\sqrt{\pi}}{2^{\frac{1}{2}\nu-\frac{1}{4}}\Gamma\left(\frac{1}{4}+\frac{1}{2}\nu\right)}\period
%\end{align}

%% References


\bibliographystyle{elsarticle-num}
%\bibliography{mybibfile}



\end{document}